\documentclass[11pt]{article}
\usepackage{epsfig}

\def\be{\begin{equation}}
\def\ee{\end{equation}}
\def\bea{\begin{eqnarray}}
\def\eea{\end{eqnarray}}
\def\bes{\begin{eqnarray*}}
\def\ees{\end{eqnarray*}}

\def\nn{\nonumber}
\def\<{\langle}
\def\>{\rangle}
\def\lb{\label}
\def\bs{\setminus}
\def\pt{\partial}

\def\R{{\bf R}}
\def\C{{\bf C}}
\def\Z{{\bf Z}}
\def\N{{\bf N}}
\def\U{{\bf U}}
\def\Q{{\bf Q}}

\def\ga{{\gamma}}

\def\ka{{\kappa}}
\def\th{{\theta}}
\def\om{{\omega}}
\def\Om{{\Omega}}
\def\ep{{\epsilon}}
\def\lm{{\lambda}}
\def\Lm{{\Lambda}}

\def\sg{{\sigma}}

\def\P{{\cal P}}

\def\rank{{\rm rank}}
\def\Sp{{\rm Sp}}
\def\mod{{\rm mod}}

\def\ol#1{\overline{#1}}  

\def\hb{\vrule height0.18cm width0.14cm $\,$}

\def\ol#1{\overline{#1}}  

\title{Multiple closed geodesics on  positively curved Finsler manifolds}
\author{
Wei Wang\thanks{Partially supported by NSFC No. 11222105, 11431001, E-mail: alexanderweiwang@gmail.com}\\\\
 School of Mathematical Sciences and LMAM \\Peking University, Beijing 100871\\
The People's Republic of China\\}
\date{}
\begin{document}

\maketitle

\begin{abstract}
{\it In this paper, we prove that on every Finsler  manifold
$(M,\,F)$  with reversibility $\lambda$ and flag
curvature $K$ satisfying
$\left(\frac{\lambda}{\lambda+1}\right)^2<K\le 1$,  there exist
$[\frac{\dim M+1}{2}]$ closed geodesics. If  the number of closed geodesics is finite,
then there exist $[\frac{\dim M}{2}]$ non-hyperbolic closed geodesics. 
Moreover, there are 3 closed geodesics on $(M,\,F)$ satisfying the above pinching condition 
when $\dim M=3$.}
\end{abstract}



{\bf Running head}: Closed geodesics on Finsler manifolds

\renewcommand{\theequation}{\thesection.\arabic{equation}}
\renewcommand{\thefigure}{\thesection.\arabic{figure}}

\setcounter{equation}{0}
\section{Introduction and main results}

This paper is devoted to a study on closed geodesics on Finsler
manifolds.  Let us recall firstly the definition of the Finsler
metrics.

{\bf Definition 1.1.} (cf. \cite{She1}) {\it Let $M$ be a finite
dimensional manifold. A function $F:TM\to [0,+\infty)$ is a {\rm
Finsler metric} if it satisfies

(F1) $F$ is $C^{\infty}$ on $TM\bs\{0\}$,

(F2) $F(x,\lm y) = \lm F(x,y)$ for all $y\in T_xM$, $x\in M$, and
$\lm>0$,

(F3) For every $y\in T_xM\bs\{0\}$, the quadratic form
$$ g_{x,y}(u,v) \equiv
         \frac{1}{2}\frac{\pt^2}{\pt s\pt t}F^2(x,y+su+tv)|_{t=s=0},
         \qquad \forall u, v\in T_xM, $$
is positive definite.

In this case, $(M,F)$ is called a {\rm Finsler manifold}. $F$ is
{\rm reversible} if $F(x,-y)=F(x,y)$ holds for all $y\in T_xM$ and
$x\in M$. $F$ is {\rm Riemannian} if $F(x,y)^2=\frac{1}{2}G(x)y\cdot
y$ for some symmetric positive definite matrix function $G(x)\in
GL(T_xM)$ depending on $x\in M$ smoothly. }

A closed curve in a Finsler manifold is a closed geodesic if it is
locally the shortest path connecting any two nearby points on this
curve (cf. \cite{She1}). As usual, on any Finsler manifold
$M=(M, F)$, a closed geodesic $c:S^1=\R/\Z\to M$ is {\it prime}
if it is not a multiple covering (i.e., iteration) of any other
closed geodesics. Here the $m$-th iteration $c^m$ of $c$ is defined
by $c^m(t)=c(mt)$. The inverse curve $c^{-1}$ of $c$ is defined by
$c^{-1}(t)=c(1-t)$ for $t\in \R$. We call two prime closed geodesics
$c$ and $d$ {\it distinct} if there is no $\th\in (0,1)$ such that
$c(t)=d(t+\th)$ for all $t\in\R$. We shall omit the word {\it
distinct} when we talk about more than one prime closed geodesic.
On a symmetric Finsler (or Riemannian) manifold, two closed geodesics
$c$ and $d$ are called { \it geometrically distinct} if $
c(S^1)\neq d(S^1)$, i.e., their image sets in $M$ are distinct.

For a closed geodesic $c$ on $(M,\,F)$, denote by $P_c$
the linearized Poincar\'{e} map of $c$  (cf. p.143 of \cite{Zil1}).
Then $P_c\in \Sp(2n-2)$ is a symplectic matrix.
For any $M\in \Sp(2k)$, we define the {\it elliptic height } $e(M)$
of $M$ to be the total algebraic multiplicity of all eigenvalues of
$M$ on the unit circle $\U=\{z\in\C|\; |z|=1\}$ in the complex plane
$\C$. Since $M$ is symplectic, $e(M)$ is even and $0\le e(M)\le 2k$.
Then $c$ is called {\it hyperbolic} if all the eigenvalues of $P_c$ avoid
$\U$, i.e., $e(P_c)=0$; {\it elliptic} if all the eigenvalues of
$P_c$ are on $\U$, i.e., $e(P_c)=2(n-1)$.

Following H.-B. Rademacher in
\cite{Rad4}, the reversibility $\lambda=\lambda(M,\,F)$ of a compact
Finsler manifold $(M,\,F)$ is defied to be
$$\lambda:=\max\{F(-X)\,|\,X\in TM, \,F(X)=1\}\ge 1.$$

It was quite surprising when A. Katok \cite{Kat1} in 1973 found
some non-symmetric Finsler metrics on CROSSs (compact rank one symmetric spaces)  with only finitely
many prime closed geodesics and all closed geodesics are
non-degenerate and elliptic. In Katok's examples the spheres $S^{2n}$ and
$S^{2n-1}$ have precisely $2n$ closed geodesics (cf. also \cite{Zil1}).

 We are aware of a number of results concerning closed geodesics on Finsler manifolds. According to the classical theorem of Lyusternik-Fet
\cite{LyF} from 1951, there exists at least one closed geodesic on
every compact Riemannian manifold. The proof of this theorem is
variational and carries over to the Finsler case.  In \cite{BaL}, V. Bangert and Y. Long proved that on any
Finsler 2-sphere $(S^2, F)$, there exist at least two  closed
geodesics. In \cite{Rad5}, H.-B. Rademacher studied the existence and
stability of closed geodesics on positively curved Finsler
manifolds. In \cite{Wan1}-\cite{Wan3}, W. Wang studied the existence and
stability of closed geodesics on positively curved Finsler spheres. 
In \cite{DuL1} of Duan and Long and in  \cite{Rad6} of Rademacher, they proved there exist 
at least two  closed geodesics on any bumpy
Finsler $n$-sphere independently. 
In \cite{LoD} and \cite{DuL2} of Duan and Long, they proved there exist 
at least two  closed geodesics on any compact simply-connected
Finsler $3$ and $4$ manifold. In  \cite{Rad7}, Rademacher proved there exist at least two  closed geodesics on any bumpy Finsler $\bf{CP^2}$.
In \cite{DLW1}, Duan, Long and Wang proved there exist at least two 
closed geodesics on any compact simply-connected bumpy Finsler manifold. 
In \cite{DLW2}, Duan, Long and Wang proved there exist at least $\frac{dn(n+1)}{2}$ 
 non-hyperbolic closed geodesics on any compact simply-connected bunpy Finsler manifold $(M,F)$  with $H^*(M;\Q)\cong T_{d,n+1}(x)$ and $K\ge 0$. In \cite{LiX}, Liu and Xiao proved there exist at least two non-contractible  closed geodesics on
any bumpy Finsler ${\rm RP}^n$.
In \cite{LLX}, Liu, Long and Xiao
proved there exist at least two non-contractible  closed geodesics on any bumpy Finsler $S^n/\Gamma$,
where $\Gamma$ is a finite group acts on $S^n$ freely  and isometrically.
In \cite{Liu}, Liu proved there exist $2[\frac{n+1}{2}]$  closed geodesics on any bumpy Finsler $S^n/\Gamma$ under a pinching condition.  In \cite{GG} and \cite{GGM} V. Ginzburg-B. Gurel and V. Ginzburg-B. Gurel-L. Macarini obtains some lower bounds for the number of closed geodesics
on certain Finsler manifolds under some conditions on the initial index of closed geodesics.

The following are the main results in this paper:

{\bf Theorem 1.2.} {\it On every Finsler  manifold
$(M,\,F)$  with reversibility $\lambda$ and flag
curvature $K$ satisfying
$\left(\frac{\lambda}{\lambda+1}\right)^2<K\le 1$,  there exist
$[\frac{\dim M+1}{2}]$ closed geodesics. }

{\bf Theorem 1.3.} {\it On every  Finsler  manifold
$(M,\,F)$  with reversibility $\lambda$ and flag
curvature $K$ satisfying
$\left(\frac{\lambda}{\lambda+1}\right)^2<K\le 1$,  there exist
$[\frac{\dim M}{2}]$ non-hyperbolic closed geodesics provided 
the number of closed geodesics on $(M, F)$ is finite. }

{\bf Theorem 1.4.} {\it On every  Finsler  manifold
$(M,\,F)$  with reversibility $\lambda$ and flag
curvature $K$ satisfying
$\left(\frac{\lambda}{\lambda+1}\right)^2<K\le 1$,  there exist
$[\frac{\dim M-1}{2}]$ closed geodesics whose linearized Poincar\'e map  
possess an eigenvalue of the form $\exp(\pi i \mu)$ with an
irrational $\mu$
provided 
the number of closed geodesics on $(M, F)$ is finite. }

{\bf Theorem 1.5.} {\it On every Finsler  $3$-manifold
$(M,\,F)$ with reversibility $\lambda$ and flag
curvature $K$ satisfying
$\left(\frac{\lambda}{\lambda+1}\right)^2<K\le 1$,  there exist
$3$ closed geodesics. }

{\bf Remark 1.6.} Note that on the standard Riemannian $n$-sphere of constant
curvature $1$, all geodesics are closed and their
linearized Poincar\'e map  are $I_{2n-2}$, therefore it
possess no eigenvalue of the form $\exp(\pi i \mu)$ with an
irrational $\mu$. Thus one can not hope that
Theorems 1.4 hold for all Finsler manifolds.
Note also that in \cite{LoW} of Y. Long and the author, they
proved the existence of two closed geodesics whose linearized Poincar\'e map  
possess an eigenvalue of the form $\exp(\pi i \mu)$ with an
irrational $\mu$
 on every Finsler  $2$-sphere $(S^2,\,F)$
provided the number of  closed geodesics is finite
by a  different method.

The proof of these theorems is motivated by 
\cite{LoZ}. In this paper, we use the Fadell-Rabinowitz index
theory in a relative version to obtain the desired critical
values of the energy functional $E$ on the space pair $(\Lambda,\,\Lambda^0)$,
where $\Lambda$ is the free loop space of $M$ and $\Lambda^0$
is its subspace consisting of constant point curves.
Then we use the method of index iteration theory of
sympletic paths developed by Y. Long and his coworkers,
especially the common index jump theorem to obtain the desired results.

In this paper, let $\N$, $\N_0$, $\Z$, $\Q$, $\R$, and $\C$ denote
the sets of natural integers, non-negative integers, integers,
rational numbers, real numbers, and complex numbers respectively.
We use only singular homology modules with $\Q$-coefficients.
For an $S^1$-space $X$, we denote by $\overline{X}$ the quotient space $X/S^1$.
We define the functions
\be \left\{\matrix{[a]=\max\{k\in\Z\,|\,k\le a\}, &
           E(a)=\min\{k\in\Z\,|\,k\ge a\} , \cr
    \varphi(a)=E(a)-[a],   \cr}\right. \lb{1.1}\ee
Especially, $\varphi(a)=0$ if $ a\in\Z\,$, and $\varphi(a)=1$ if $
a\notin\Z\,$.

\setcounter{equation}{0}
\section{Critical point theory for closed geodesics}

In this section, we describe briefly the critical point theory for closed geodesics.

On a compact Finsler manifold $(M,F)$, we choose an auxiliary Riemannian
metric. This endows the space $\Lambda=\Lambda M$ of $H^1$-maps
$\gamma:S^1\rightarrow M$ with a natural Riemannian Hilbert manifold structure
on which the group $S^1=\R/\Z$ acts continuously
by isometries, cf. \cite{Kli2}, Chapters 1 and 2. This action is
defined by translating the parameter, i.e.,
$$ (s\cdot\gamma)(t)=\gamma(t+s) \qquad $$
for all $\gamma\in\Lm$ and $s,t\in S^1$.
The Finsler metric $F$ defines an energy functional $E$ and a length
functional $L$ on $\Lambda$ by
\be E(\gamma)=\frac{1}{2}\int_{S^1}F(\dot{\gamma}(t))^2dt,
 \quad L(\gamma) = \int_{S^1}F(\dot{\gamma}(t))dt.  \lb{2.1}\ee
Both functionals are invariant under the $S^1$-action.
By \cite{Mer1}, the functional $E$ is $C^{1, 1}$ on $\Lm$ and
satisfies the Palais-Smale condition. Thus we can apply the
deformation theorems in \cite{Cha} and \cite{MaW}.
The critical points
of $E$ of positive energies are precisely the closed geodesics $c:S^1\to M$
of the Finsler structure. If $c\in\Lambda$ is a closed geodesic then $c$ is
a regular curve, i.e., $\dot{c}(t)\not= 0$ for all $t\in S^1$, and this implies
that the second differential $E''(c)$ of $E$ at $c$ exists.
As usual we define the index $i(c)$ of $c$ as the maximal dimension of
subspaces of $T_c \Lambda$ on which $E^{\prime\prime}(c)$ is negative definite, and the
nullity $\nu(c)$ of $c$ so that $\nu(c)+1$ is the dimension of the null
space of $E^{\prime\prime}(c)$.

For $m\in\N$ we denote the $m$-fold iteration map
$\phi^m:\Lambda\rightarrow\Lambda$ by \be \phi^m(\ga)(t)=\ga(mt)
\qquad \forall\,\ga\in\Lm, t\in S^1. \lb{2.2}\ee We also use the
notation $\phi^m(\gamma)=\gamma^m$. For a closed geodesic $c$, the
average index is defined by \be
\hat{i}(c)=\lim_{m\rightarrow\infty}\frac{i(c^m)}{m}. \lb{2.3}\ee

If $\gamma\in\Lambda$ is not constant then the multiplicity
$m(\gamma)$ of $\gamma$ is the order of the isotropy group $\{s\in
S^1\mid s\cdot\gamma=\gamma\}$. If $m(\gamma)=1$ then $\gamma$ is
called {\it prime}. Hence $m(\gamma)=m$ if and only if there exists a
prime curve $\tilde{\gamma}\in\Lambda$ such that
$\gamma=\tilde{\gamma}^m$.

In this paper for $\ka\in \R$ we denote by
\be \Lm^{\ka}=\{d\in \Lm\,|\,E(d)\le \ka\}.   \lb{2.4}\ee
For a closed geodesic $c$ we set
$$ \Lm(c)=\{\ga\in\Lm\mid E(\ga)<E(c)\}. $$

We call a closed geodesic satisfying the isolation condition, if
the following holds:

{\bf (Iso)  For all $m\in\N$ the orbit $S^1\cdot c^m$ is an
isolated critical orbit of $E$. }

Note that if the number of prime closed geodesics on a Finsler manifold
is finite, then all the closed geodesics satisfy (Iso).

Using singular homology with rational
coefficients we consider the following critical $\Q$-module of a closed geodesic
$c\in\Lambda$:
\be \overline{C}_*(E,c)
   = H_*\left((\Lm(c)\cup S^1\cdot c)/S^1,\Lm(c)/S^1\right). \lb{2.5}\ee

{\bf Proposition 2.1.} (cf. Satz 6.11 of \cite{Rad2} or Proposition
3.12 of \cite{BaL}) {\it Let $c$ be a prime closed geodesic on a Finsler
manifold $(M,F)$ satisfying (Iso). Then we have
\bea \overline{C}_q( E,c^m)
&&\equiv H_q\left((\Lm(c^m)\cup S^1\cdot c^m)/S^1, \Lm(c^m)/S^1\right)\nn\\
&&= \left(H_{i(c^m)}(U_{c^m}^-\cup\{c^m\},U_{c^m}^-)
    \otimes H_{q-i(c^m)}(N_{c^m}^-\cup\{c^m\},N_{c^m}^-)\right)^{+\Z_m} \nn
\eea

(i) When $\nu(c^m)=0$, there holds
\bea \overline{C}_q( E,c^m) = \left\{\matrix{
     \Q, &\quad {\it if}\;\; i(c^m)-i(c)\in 2\Z,\;{\it and}\;
                   q=i(c^m),\;  \cr
     0, &\quad {\it otherwise}. \cr}\right.  \nn\eea

(ii) When $\nu(c^m)>0$, there holds
\bea \overline{C}_q( E,c^m) =
H_{q-i(c^m)}(N_{c^m}^-\cup\{c^m\}, N_{c^m}^-)^{(-1)^{i(c^m)-i(c)}\Z_m},\nn\eea
where $N_{c^m}$ is a local characteristic manifold at $c^m$ and $N^-_{c^m}=N_{c^m}\cap \Lambda(c^m)$,  $U_{c^m}$ 
is a local negative disk at $c^m$ and $U^-_{c^m}=U_{c^m}\cap \Lambda(c^m)$,
$H_{\ast}(X,A)^{\pm\Z_m}
   = \{[\xi]\in H_{\ast}(X,A)\,|\,T_{\ast}[\xi]=\pm [\xi]\}$
where $T$ is a generator of the $\Z_m$-action.}

Denote by 
\bea k_j(c^m) \equiv \dim\, H_j(N_{c^m}^-\cup\{c^m\},N_{c^m}^- )^{(-1)^{i(c^m)-i(c)}\Z_m}. \lb{2.6}
\eea

Clearly the integers $k_j(c^m)$  equal to
$0$ when $j<0$ or $j>\nu(c^m)$ and can take only values $0$ or $1$
when $j=0$ or $j=\nu(c^m)$.

{\bf Proposition 2.2.} (cf. Satz 6.13 of \cite{Rad2}) {\it Let $c$
be a prime closed geodesic on a Finsler manifold $(M, F)$
satisfying (Iso). For any $m\in\N$, we have

(i) If $k_0(c^m)=1$, there holds $k_j(c^m)=0$ for
$1\le j\le \nu(c^m)$.

(ii) If $k_{\nu(c^m)}(c^m)=1$,
there holds $k_j(c^m)=0$ for $0\le j\le \nu(c^m)-1$.

(iii) If $k_j(c^m)\ge 1$  for some
$1\le j\le \nu(c^m)-1$, there holds $k_{\nu(c^m)}(c^m)=0=k_0(c^m)$.

(iv) In particular, if $\nu(c^m)\le 2$, then only one of the
$k_j(c^m)$'s can be non-zero.}

By Lemma 5.2 of \cite{Wan1}, we have  the following  periodic property for $k_l(c^m)$:

{\bf Proposition 2.3.} {\it  Let $c$ be a prime closed geodesic on a
compact Finsler manifold $(M, F)$ satisfying (Iso).
Then there exists a minimal $T(c)\in \N$ such that
\bea
&& \nu(c^{p+T(c)})=\nu(c^p),\quad i(c^{p+T(c)})-i(c^p)\in 2\Z,
     \qquad\forall p\in \N,  \lb{2.7}\\
&& k_l(c^{p+T(c)})=k_l(c^p), \qquad\forall p\in \N,\;l\in\Z. \lb{2.8}\eea
}

{\bf Definition 2.4.} {\it The Euler characteristic $\chi(c^m)$
of $c^m$ is defined by
\bea \chi(c^m)
&\equiv& \chi\left((\Lm(c^m)\cup S^1\cdot c^m)/S^1, \Lm(c^m)/S^1\right), \nn\\
&\equiv& \sum_{q=0}^{\infty}(-1)^q\dim \overline{C}_q( E,c^m)
=\sum_{l=0}^{2n-2}(-1)^{i(c^m)+l}k_l(c^m).\lb{2.9}\eea
Here $\chi(A, B)$ denotes the usual Euler characteristic of the space pair $(A, B)$.

The average Euler characteristic $\hat\chi(c)$ of $c$ is defined by }
\be \hat{\chi}(c)=\lim_{N\to\infty}\frac{1}{N}\sum_{1\le m\le N}\chi(c^m).
\lb{2.10}\ee
By Remark 5.4 of \cite{Wan1},  $\hat\chi(c)$ is well-defined and is a
rational number, in fact \bea
\hat\chi(c) 
=\frac{1}{T(c)}\sum_{1\le m\le T(c)}\chi(c^m).\lb{2.11}\eea
In particular, if $c^m$ are non-degenerate for $\forall m\in\N$, then 
 \bea
 \hat\chi(c)=\;\;\left\{\matrix{
    (-1)^{i(c)},&\quad {\it if}\quad i(c^2)-i(c)\in2\Z,  \cr
    \frac{(-1)^{i(c)}}{2},&\quad {\it otherwise.}\cr
    }\right.  \lb{2.12}\eea

We have the following mean index identity for closed geodesics:

{\bf Theorem 2.5.}  (cf. Theorem 7.9 in \cite{Rad2} or Theorem 5.5 in \cite{Wan1}) {\it Suppose that there exist only finitely many prime
closed geodesics $\{c_j\}_{1\le j\le p}$ with $\hat i(c_j)>0$
for $1\le j\le p$ on $(S^n,F)$.
Then the following identity holds
\bea  \sum_{1\le j\le p}\frac{\hat\chi(c_j)}{\hat{i}(c_j)}=B(n, 1)=\left\{\matrix{
     \frac{-n}{2n-2}, &\quad n\quad even, \cr
     \frac{n+1}{2n-2} &\quad n\quad odd. \cr}\right.
     \lb{2.13}\eea
}

Set $\ol{\Lm}^0=\ol{\Lambda}^0M =\{{\rm
constant\;point\;curves\;in\;}M\}\cong M$. Let $(X,Y)$ be a
space pair such that the Betti numbers $b_i=b_i(X,Y)=\dim
H_i(X,Y;\Q)$ are finite for all $i\in \Z$. As usual the {\it
Poincar\'e series} of $(X,Y)$ is defined by the formal power series
$P(X, Y)=\sum_{i=0}^{\infty}b_it^i$. We need the following  results on Betti numbers.

{\bf Theorem 2.6.} (H.-B. Rademacher, Theorem 2.4 and Remark 2.5 of \cite{Rad1})
{\it We have the Poincar\'e series

(i) When $n=2k+1$ is odd
\bea
P(\ol{\Lm}S^n,\ol{\Lm}^0S^n)(t)
&=&t^{n-1}\left(\frac{1}{1-t^2}+\frac{t^{n-1}}{1-t^{n-1}}\right)\nn\\
&=& t^{2k}\left(\frac{1}{1-t^2}+\frac{t^{2k}}{1- t^{2k}}\right).
\lb{2.14}\eea
Thus for $q\in\Z$ and $l\in\N_0$, we have
\bea {b}_q &=& {b}_q(\ol{\Lm}S^n,\ol{\Lm}^0 S^n)\nn\\
 &=&\rank H_q(\ol{\Lm} S^n,\ol{\Lm}^0 S^n )\nn\\
 &=& \;\;\left\{\matrix{
    2,&\quad {\it if}\quad q\in \{4k+2l,\quad l=0\;\mod\; k\},  \cr
    1,&\quad {\it if}\quad q\in \{2k\}\cup\{2k+2l,\quad l\neq 0\;\mod\; k\},  \cr
    0 &\quad {\it otherwise}. \cr}\right. \lb{2.15}\eea

(ii) When $n=2k$ is even
\bea
P(\ol{\Lm}S^n,\ol{\Lm}^0S^n)(t)
&=&t^{n-1}\left(\frac{1}{1-t^2}+\frac{t^{n(m+1)-2}}{1-t^{n(m+1)-2}}\right)
\frac{1-t^{nm}}{1-t^n}\nn\\
&=& t^{2k-1}\left(\frac{1}{1-t^2}+\frac{t^{4k-2}}{1- t^{4k-2}}\right),
\lb{2.16}\eea
where $m=1$ by Theorem 2.4 of \cite{Rad1}. Thus for $q\in\Z$ and $l\in\N_0$, we have}
\bea {b}_q &=& {b}_q(\ol{\Lm}S^n,\ol{\Lm}^0 S^n)\nn\\
 &=&\rank H_q(\ol{\Lm} S^n,\ol{\Lm}^0 S^n )\nn\\
 &=& \;\;\left\{\matrix{
    2,&\quad {\it if}\quad q\in \{6k-3+2l,\quad l=0\;\mod\; 2k-1\},  \cr
    1,&\quad {\it if}\quad q\in \{2k-1\}\cup\{2k-1+2l,\quad l\neq 0\;\mod\; 2k-1\},  \cr
    0 &\quad {\it otherwise}. \cr}\right. \lb{2.17}\eea

We have the following version of the Morse inequality.

{\bf Theorem 2.7.} (Theorem 6.1 of \cite{Rad2}) {\it Suppose that there exist
only finitely many prime closed geodesics $\{c_j\}_{1\le j\le p}$ on $(M, F)$,
and $0\le a<b\le \infty$.
Define for each $q\in\Z$,
\bea
{M}_q(\ol{\Lm}^b,\ol{\Lm}^a)
&=& \sum_{1\le j\le p,\;a<E(c^m_j)<b}\rank{\ol{C}}_q(E, c^m_j ) \nn\\
{b}_q(\ol{\Lm}^{b},\ol{\Lm}^{a})
&=& \rank H_q(\ol{\Lm}^{b},\ol{\Lm}^{a}). \nn\eea
Then there holds }
\bea
M_q(\ol{\Lm}^{b},\ol{\Lm}^{a}) &-& M_{q-1}(\ol{\Lm}^{b},\ol{\Lm}^{a})
    + \cdots +(-1)^{q}M_0(\ol{\Lm}^{b},\ol{\Lm}^{a}) \nn\\
&\ge& b_q(\ol{\Lm}^{b},\ol{\Lm}^{a}) - b_{q-1}(\ol{\Lm}^{b},\ol{\Lm}^{a})
   + \cdots + (-1)^{q}b_0(\ol{\Lm}^{b},\ol{\Lm}^{a}), \lb{2.18}\\
{M}_q(\ol{\Lm}^{b},\ol{\Lm}^{a}) &\ge& {b}_q(\ol{\Lm}^{b},\ol{\Lm}^{a}).\lb{2.19}
\eea

Next we recall the Fadell-Rabinowitz index
in a relative version due to \cite{Rad3}.
Let $X$ be an $S^1$-space, $A\subset X$ a closed $S^1$-invariant subset.
Note that the cup product defines a homomorphism
\be H^\ast_{S^1}(X)\otimes  H^\ast_{S^1}(X,\; A)\rightarrow  H^\ast_{S^1}(X,\;A):
\quad (\zeta,\;z)\rightarrow \zeta\cup z,\lb{2.20}\ee
where $H_{S^1}^\ast$ is the $S^1$-equivariant cohomology with
rational coefficients in the sense of A. Borel
(cf. Chapter IV of \cite{Bor1}).
We fix a characteristic class $\eta\in H^2(CP^\infty)$. Let
$f^\ast: H^\ast(CP^\infty)\rightarrow H_{S^1}^\ast(X)$ be the
homomorphism induced by a classifying map $f: X_{S^1}\rightarrow CP^\infty$.
Now for $\gamma\in  H^\ast(CP^\infty)$ and $z\in H^\ast_{S^1}(X,\;A)$, let
$\gamma\cdot z=f^\ast(\gamma)\cup z$. Then
the order $ord_\eta(z)$ with respect to $\eta$ is defined by
\be ord_\eta(z)=\inf\{k\in\N\cup\{\infty\}\;|\;\eta^k\cdot z= 0\}.\lb{2.21}\ee
By Proposition 3.1 of \cite{Rad3}, there is an element
$z\in H^{n+1}_{S^1}(\Lambda,\;\Lambda^0)$ of infinite order, i.e., $ord_\eta(z)=\infty$.
For $\kappa\ge0$, we denote by
$j_\kappa: (\Lambda^\kappa, \;\Lambda^0)\rightarrow (\Lambda \;\Lambda^0)$
the natural inclusion and define the function $d_z:\R^{\ge 0}\rightarrow \N\cup\{\infty\}$:
\be d_z(\kappa)=ord_\eta(j_\kappa^\ast(z)).\lb{2.22}\ee
Denote by $d_z(\kappa-)=\lim_{\epsilon\searrow 0}d_z(\kappa-\epsilon)$,
where  $t\searrow a$  means $t>a$ and $t\to a$.

Then we have the following property due to Section 5 of \cite{Rad3}:

{\bf Lemma 2.8.} (H.-B. Rademacher) {\it The function $d_z$ is non-decreasing and
$\lim_{\lambda\searrow \kappa}d_z(\lambda)=d_z(\kappa)$.
 Each discontinuous point of $d_z$ is a critical value of the energy functional $E$.
In particular, if $d_z(\kappa)-d_z(\kappa-)\ge 2$, then there are infinitely
many prime closed geodesics $c$ with energy $\kappa$.
}\hfill\hb

For each $i\ge 1$, we define
\be \kappa_i=\inf\{\delta\in\R\;|\: d_z(\delta)\ge i\}.\lb{2.23}\ee
Then we have the following:

{\bf Lemma 2.9.}  (cf. Lemma 2.3 of \cite{Wan3}) {\it Suppose there are only finitely many prime closed
geodesics on $(S^n,\, F)$. Then each $\kappa_i$ is a critical value of $E$.
If $\kappa_i=\kappa_j$ for some $i<j$, then there are infinitely many prime
closed geodesics on $(S^n,\, F)$. }

{\bf Lemma 2.10.}  (cf. Lemma 2.4 of \cite{Wan3}) {\it Suppose there are only finitely many prime closed
geodesics on $(S^n,\, F)$. Then for every $i\in\N$, there exists a
closed geodesic $c$ on $(S^n,\, F)$ such that
\bea E(c)=\kappa_i,\quad
\ol{C}_{2i+\dim(z)-2}(E, c)\neq 0.\lb{2.24}\eea}

{\bf Definition 2.11.} {\it A prime closed geodesic $c$ is
$(m, i)$-{ \bf variationally visible:} if there exist some $m, i\in\N$
such that (\ref{2.24}) holds for $c^m$ and $\kappa_i$. We call $c$
{\bf infinitely variationally visible:} if there exist
infinitely many $m, i\in\N$ such that $c$ is $(m, i)$-variationally visible.
We denote by $\mathcal{V}_\infty(S^n, F)$ the set of infinitely
variationally visible closed geodesics.}

{\bf Theorem 2.12.}  (cf. Theorem 2.6 of \cite{Wan3}) {\it Suppose there are only finitely many prime closed
geodesics on $(S^n,\, F)$. Then for any $c\in\mathcal{V}_\infty(S^n, F)$,
we have
\be \frac{\hat i(c)}{L(c)}=2\sigma.\lb{2.25}
\ee
where $\sigma=\lim\inf_{i\rightarrow\infty}i/\sqrt{2\kappa_i}=
\lim\sup_{i\rightarrow\infty}{i}/{\sqrt{2\kappa_i}}$.
}

\setcounter{equation}{0}
\section{ Index iteration theory for closed geodesics}

In this section, we recall briefly the index theory for symplectic paths
developed by Y. Long and his coworkers. All the details can be found in \cite{Lon4}.
Then we use this theory to study the Morse indices of closed geodesics.

Let $c$ be a closed geodesic on an orientable Finsler manifold $M=(M,\,F)$.
Denote the linearized Poincar\'e map of $c$ by $P_c\in\Sp(2n-2)$.
Then $P_c$ is a symplectic matrix.
Note that the index iteration formulae in \cite{Lon3} of 2000 (cf. Chap. 8 of
\cite{Lon4}) work for Morse indices of iterated closed geodesics (cf.
\cite{LLo}, Chap. 12 of \cite{Lon4}). Since every closed geodesic
on $M$ is orientable. Then by Theorem 1.1 of \cite{Liu1}
of C. Liu (cf. also \cite{Wil}), the initial Morse index of a closed geodesic
$c$ on $M$ coincides with the index of a
corresponding symplectic path introduced by C. Conley, E. Zehnder, and Y. Long
in 1984-1990 (cf. \cite{Lon4}).

As usual, the symplectic group $\Sp(2n)$ is defined by
$$ \Sp(2n) = \{M\in {\rm GL}(2n,\R)\,|\,M^TJM=J\}, $$
whose topology is induced from that of $\R^{4n^2}$. For $\tau>0$ we are interested
in paths in $\Sp(2n)$:
$$ \P_{\tau}(2n) = \{\ga\in C([0,\tau],\Sp(2n))\,|\,\ga(0)=I_{2n}\}, $$
which is equipped with the topology induced from that of $\Sp(2n)$. The
following real function was introduced in \cite{Lon3}:
$$ D_{\om}(M) = (-1)^{n-1}\ol{\om}^n\det(M-\om I_{2n}), \qquad
          \forall \om\in\U,\, M\in\Sp(2n). $$
Thus for any $\om\in\U$ the following codimension $1$ hypersurface in $\Sp(2n)$ is
defined in \cite{Lon3}:
$$ \Sp(2n)_{\om}^0 = \{M\in\Sp(2n)\,|\, D_{\om}(M)=0\}.  $$
For any $M\in \Sp(2n)_{\om}^0$, we define a co-orientation of $\Sp(2n)_{\om}^0$
at $M$ by the positive direction $\frac{d}{dt}Me^{t\ep J}|_{t=0}$ of
the path $Me^{t\ep J}$ with $0\le t\le 1$ and $\ep>0$ being sufficiently
small. Let
\bea
\Sp(2n)_{\om}^{\ast} &=& \Sp(2n)\bs \Sp(2n)_{\om}^0,   \nn\\
\P_{\tau,\om}^{\ast}(2n) &=&
      \{\ga\in\P_{\tau}(2n)\,|\,\ga(\tau)\in\Sp(2n)_{\om}^{\ast}\}, \nn\\
\P_{\tau,\om}^0(2n) &=& \P_{\tau}(2n)\bs  \P_{\tau,\om}^{\ast}(2n).  \nn\eea
For any two continuous arcs $\xi$ and $\eta:[0,\tau]\to\Sp(2n)$ with
$\xi(\tau)=\eta(0)$, it is defined as usual:
$$ \eta\ast\xi(t) = \left\{\matrix{
            \xi(2t), & \quad {\rm if}\;0\le t\le \tau/2, \cr
            \eta(2t-\tau), & \quad {\rm if}\; \tau/2\le t\le \tau. \cr}\right. $$
Given any two $2m_k\times 2m_k$ matrices of square block form
$M_k=\left(\matrix{A_k&B_k\cr
                                C_k&D_k\cr}\right)$ with $k=1, 2$,
as in \cite{Lon4}, the $\;\diamond$-product of $M_1$ and $M_2$ is defined by
the following $2(m_1+m_2)\times 2(m_1+m_2)$ matrix $M_1\diamond M_2$:
$$ M_1\diamond M_2=\left(\matrix{A_1&  0&B_1&  0\cr
                               0&A_2&  0&B_2\cr
                             C_1&  0&D_1&  0\cr
                               0&C_2&  0&D_2\cr}\right). \nn$$  
Denote by $M^{\diamond k}$ the $k$-fold $\diamond$-product $M\diamond\cdots\diamond M$. Note
that the $\diamond$-product of any two symplectic matrices is symplectic. For any two
paths $\ga_j\in\P_{\tau}(2n_j)$ with $j=0$ and $1$, let
$\ga_0\diamond\ga_1(t)= \ga_0(t)\diamond\ga_1(t)$ for all $t\in [0,\tau]$.

A special path $\xi_n$ is defined by
\be \xi_n(t) = \left(\matrix{2-\frac{t}{\tau} & 0 \cr
                                             0 &  (2-\frac{t}{\tau})^{-1}\cr}\right)^{\diamond n}
        \qquad {\rm for}\;0\le t\le \tau.  \lb{3.1}\ee

{\bf Definition 3.1.} (cf. \cite{Lon3}, \cite{Lon4}) {\it For any $\om\in\U$ and
$M\in \Sp(2n)$, define
\be  \nu_{\om}(M)=\dim_{\C}\ker_{\C}(M - \om I_{2n}).  \lb{3.2}\ee
For any $\tau>0$ and $\ga\in \P_{\tau}(2n)$, define
\be  \nu_{\om}(\ga)= \nu_{\om}(\ga(\tau)).  \lb{3.3}\ee

If $\ga\in\P_{\tau,\om}^{\ast}(2n)$, define
\be i_{\om}(\ga) = [\Sp(2n)_{\om}^0: \ga\ast\xi_n],  \lb{3.4}\ee
where the right hand side of (\ref{3.4}) is the usual homotopy intersection
number, and the orientation of $\ga\ast\xi_n$ is its positive time direction under
homotopy with fixed end points.

If $\ga\in\P_{\tau,\om}^0(2n)$, we let $\mathcal{F}(\ga)$
be the set of all open neighborhoods of $\ga$ in $\P_{\tau}(2n)$, and define
\be i_{\om}(\ga) = \sup_{U\in\mathcal{F}(\ga)}\inf\{i_{\om}(\beta)\,|\,
                       \beta\in U\cap\P_{\tau,\om}^{\ast}(2n)\}.
               \lb{3.5}\ee
Then
$$ (i_{\om}(\ga), \nu_{\om}(\ga)) \in \Z\times \{0,1,\ldots,2n\}, $$
is called the index function of $\ga$ at $\om$. }

For any symplectic path $\ga\in\P_{\tau}(2n)$ and $m\in\N$,  we
define its $m$-th iteration $\ga^m:[0,m\tau]\to\Sp(2n)$ by
\be \ga^m(t) = \ga(t-j\tau)\ga(\tau)^j, \qquad
  {\rm for}\quad j\tau\leq t\leq (j+1)\tau,\;j=0,1,\ldots,m-1.
     \lb{3.6}\ee
We still denote the extended path on $[0,+\infty)$ by $\ga$.

{\bf Definition 3.2.} (cf. \cite{Lon3}, \cite{Lon4}) {\it For any $\ga\in\P_{\tau}(2n)$,
we define
\be (i(\ga,m), \nu(\ga,m)) = (i_1(\ga^m), \nu_1(\ga^m)), \qquad \forall m\in\N.
   \lb{3.7}\ee
The mean index $\hat{i}(\ga,m)$ per $m\tau$ for $m\in\N$ is defined by
\be \hat{i}(\ga,m) = \lim_{k\to +\infty}\frac{i(\ga,mk)}{k}. \lb{3.8}\ee
For any $M\in\Sp(2n)$ and $\om\in\U$, the {\it splitting numbers} $S_M^{\pm}(\om)$
of $M$ at $\om$ are defined by
\be S_M^{\pm}(\om)
     = \lim_{\ep\to 0^+}i_{\om\exp(\pm\sqrt{-1}\ep)}(\ga) - i_{\om}(\ga),
   \lb{3.9}\ee
for any path $\ga\in\P_{\tau}(2n)$ satisfying $\ga(\tau)=M$.}

For a given path $\gamma\in {\cal P}_{\tau}(2n)$ we consider to deform
it to a new path $\eta$ in ${\cal P}_{\tau}(2n)$ so that
\begin{equation}
i_1(\gamma^m)=i_1(\eta^m),\quad \nu_1(\gamma^m)=\nu_1(\eta^m), \quad
         \forall m\in {\bf N}, \label{3.10}
\end{equation}
and that $(i_1(\eta^m),\nu_1(\eta^m))$ is easy enough to compute. This
leads to finding homotopies $\delta:[0,1]\times[0,\tau]\to {\rm Sp}(2n)$
starting from $\gamma$ in ${\cal P}_{\tau}(2n)$ and keeping the end
points of the homotopy always stay in a certain suitably chosen maximal
subset of ${\rm Sp}(2n)$ so that (\ref{3.10}) always holds. In fact,  this
set was first discovered in \cite{Lon3} as the path connected component
$\Omega^0(M)$ containing $M=\gamma(\tau)$ of the set
\begin{eqnarray}
  \Omega(M)=\{N\in{\rm Sp}(2n)\,&|&\,\sigma(N)\cap{\bf U}=\sigma(M)\cap{\bf U}\;
{\rm and}\;  \nonumber\\
 &&\qquad \nu_{\lambda}(N)=\nu_{\lambda}(M),\;\forall\,
\lambda\in\sigma(M)\cap{\bf U}\}. \label{3.11}
\end{eqnarray}
Here $\Omega^0(M)$ is called the {\it homotopy component} of $M$ in
${\rm Sp}(2n)$.

In \cite{Lon3} and \cite{Lon4}, the following symplectic matrices were introduced
as {\it basic normal forms}:
\begin{eqnarray}
D(\lambda)=\left(\matrix{\lm & 0\cr
         0  & \lm^{-1}\cr}\right), &\quad& \lm=\pm 2,\lb{3.12}\\
N_1(\lm,b) = \left(\matrix{\lm & b\cr
         0  & \lm\cr}\right), &\quad& \lm=\pm 1, b=\pm1, 0, \lb{3.13}\\
R(\th)=\left(\matrix{\cos\th & -\sin\th\cr
        \sin\th  & \cos\th\cr}\right), &\quad& \th\in (0,\pi)\cup(\pi,2\pi),
                     \lb{3.14}\\
N_2(\om,b)= \left(\matrix{R(\th) & b\cr
              0 & R(\th)\cr}\right), &\quad& \th\in (0,\pi)\cup(\pi,2\pi),
                     \lb{3.15}\end{eqnarray}
where $b=\left(\matrix{b_1 & b_2\cr
               b_3 & b_4\cr}\right)$ with  $b_i\in\R$ and  $b_2\not=b_3$.
We call $ N_2(\omega, b)$ 
{\it trivial} if $(b_2-b_3)\sin\theta>0$,   $ N_2(\omega, b)$ 
{\it non-trivial} if $(b_2-b_3)\sin\theta<0$.

Splitting numbers possess the following properties:

{\bf Lemma 3.3.} (cf. \cite{Lon3} and Lemma 9.1.5 of \cite{Lon4}) {\it Splitting
numbers $S_M^{\pm}(\om)$ are well defined, i.e., they are independent of the choice
of the path $\ga\in\P_\tau(2n)$ satisfying $\ga(\tau)=M$ appeared in (\ref{3.9}).
For $\om\in\U$ and $M\in\Sp(2n)$, splitting numbers $S_N^{\pm}(\om)$ are constant
for all $N\in\Om^0(M)$. }

{\bf Lemma 3.4.} (cf. \cite{Lon3}, Lemma 9.1.5 and List 9.1.12 of \cite{Lon4})
{\it For $M\in\Sp(2n)$ and $\om\in\U$, there hold
\begin{eqnarray}
S_M^{\pm}(\om) &=& 0, \qquad {\it if}\;\;\om\not\in\sg(M).  \lb{3.16}\\
S_{N_1(1,a)}^+(1) &=& \left\{\matrix{1, &\quad {\rm if}\;\; a\ge 0, \cr
0, &\quad {\rm if}\;\; a< 0. \cr}\right. \lb{3.17}\eea

For any $M_i\in\Sp(2n_i)$ with $i=0$ and $1$, there holds }
\be S^{\pm}_{M_0\diamond M_1}(\om) = S^{\pm}_{M_0}(\om) + S^{\pm}_{M_1}(\om),
    \qquad \forall\;\om\in\U. \lb{3.18}\ee

The following is the precise index iteration formulae for symplectic paths, which is due to Y. Long (cf. Theorem 8.3.1 and Corollary 8.3.2 of \cite{Lon4}).

{\bf Theorem 3.5.} {\it Let $\gamma\in\P_\tau(2n)$, then there exists a path $f\in
C([0,1],\Omega^0(\gamma(\tau))$ such that $f(0)=\gamma(\tau)$ and
\bea f(1)=&&N_1(1,1)^{\diamond p_-} \diamond I_{2p_0}\diamond
N_1(1,-1)^{\diamond p_+}
\diamond N_1(-1,1)^{\diamond q_-} \diamond (-I_{2q_0})\diamond
N_1(-1,-1)^{\diamond q_+}\nn\\
&&\diamond R(\theta_1)\diamond\cdots\diamond R(\theta_r)
\diamond N_2(\omega_1, u_1)\diamond\cdots\diamond N_2(\omega_{r_*}, u_{r_*}) \nn\\
&&\diamond N_2(\lm_1, v_1)\diamond\cdots\diamond N_2(\lm_{r_0}, v_{r_0})
\diamond M_0 \lb{3.19}\eea
where $ N_2(\omega_j, u_j) $s are
non-trivial and   $ N_2(\lm_j, v_j)$s  are trivial basic normal
forms; $\sigma (M_0)\cap U=\emptyset$; $p_-$, $p_0$, $p_+$, $q_-$,
$q_0$, $q_+$, $r$, $r_*$ and $r_0$ are non-negative integers;
$\omega_j=e^{\sqrt{-1}\alpha_j}$, $
\lambda_j=e^{\sqrt{-1}\beta_j}$; $\theta_j$, $\alpha_j$, $\beta_j$
$\in (0, \pi)\cup (\pi, 2\pi)$; these integers and real numbers
are uniquely determined by $\gamma(\tau)$. Then using the
functions defined in (\ref{1.1})
\bea i(\gamma, m)=&&m(i(\gamma,
1)+p_-+p_0-r)+2\sum_{j=1}^r E\left(\frac{m\theta_j}{2\pi}\right)-r
-p_--p_0\nn\\&&-\frac{1+(-1)^m}{2}(q_0+q_+)+2\left(
\sum_{j=1}^{r_*}\varphi\left(\frac{m\alpha_j}{2\pi}\right)-r_*\right).
\lb{3.20}\eea
\bea \nu(\gamma, m)=&&\nu(\gamma,
1)+\frac{1+(-1)^m}{2}(q_-+2q_0+q_+)+2(r+r_*+r_0)\nn\\
&&-2\left(\sum_{j=1}^{r}\varphi\left(\frac{m\theta_j}{2\pi}\right)+
\sum_{j=1}^{r_*}\varphi\left(\frac{m\alpha_j}{2\pi}\right)
+\sum_{j=1}^{r_0}\varphi\left(\frac{m\beta_j}{2\pi}\right)\right)\lb{3.21}\eea
\bea \hat i(\gamma, 1)=i(\gamma, 1)+p_-+p_0-r+\sum_{j=1}^r
\frac{\theta_j}{\pi}.\lb{3.22}\eea
We have $i(\gamma, 1)$ is odd if $f(1)=N_1(1, 1)$, $I_2$,
$N_1(-1, 1)$, $-I_2$, $N_1(-1, -1)$ and $R(\theta)$; $i(\gamma, 1)$ is
even if $f(1)=N_1(1, -1)$ and $ N_2(\omega, b)$; $i(\gamma, 1)$ can be any
integer if $\sigma (f(1)) \cap \U=\emptyset$.}

We have the following properties in the index iteration theory.

{\bf Theorem 3.6. } (cf. Theorem 2.2 of \cite{LoZ}
)
{\it Let $\gamma\in\P_\tau(2n)$ and $M = \gamma(\tau)$. 
Then for any $m\in\N$, there holds
$$\nu(\gamma, m)-\frac{e(M)}{2}\le i(\gamma, m+1)-i(\gamma, m)-i(\gamma, 1)
\le \nu(\gamma, 1)-\nu(\gamma, m+1)+\frac{e(M)}{2},$$
where $e(M)$ is the total algebraic multiplicity of all eigenvalues of $M$ on
the unit circle in the complex plane $\C$. }

The following is the common index jump theorem of Y. Long and C. Zhu.

{\bf Theorem 3.7.} (cf. Theorems 4.1-4.3 of \cite{LoZ}
) {\it 	Let $\gamma_k\in\P_{\tau_k}(2n)$ for
$ k = 1,\ldots,q$ be a finite collection of symplectic paths. Let $M_k = \gamma_k(\tau_k)$. 
 Suppose  $\hat i(\gamma_k, 1) > 0$ for all $k = 1,\ldots,q$.
Then there exist infinitely many $(T,m_1,\ldots,m_q)\in\N^{q+1}$  such that
\bea \nu(\gamma_k, 2m_k -1)&=&	\nu(\gamma_k, 1),\lb{3.23}\\
\nu(\gamma_k, 2m_k +1)&=&	\nu(\gamma_k, 1),\lb{3.24}\\
i(\gamma_k, 2m_k -1)+\nu(\gamma_k, 2m_k -1)&=&	
2T-\left(i(\gamma_k, 1)+2S^+_{M_k}(1)-\nu(\gamma_k, 1)\right),\lb{3.25}\\
i(\gamma_k, 2m_k+1)&=&2T+i(\gamma_k, 1),\lb{3.26}\\
i(\gamma_k, 2m_k)&\ge&2T-\frac{e(M_k)}{2}\ge 2T-n,\lb{3.27}\\
i(\gamma_k, 2m_k)+\nu(\gamma_k, 2m_k)&\le&2T+\frac{e(M_k)}{2}\le 2T+n,\lb{3.28}
\eea
for every $k=1,\ldots,q$. 
Moreover we have
\bea \min\left\{\left\{\frac{m_k\theta}{\pi}\right\},\,1-\left\{\frac{m_k\theta}{\pi}\right\}\right\}
<\delta,\lb{3.29}\eea
whenever $e^{\sqrt{-1}\theta}\in\sigma(M_k)$
and $\delta$ can be chosen as small as we want (cf. (4.43) of \cite{LoZ}
). More precisely, by (4.10) and (4.40) in \cite{LoZ}
 , we have
\bea m_k=\left(\left[\frac{T}{M\hat i(\gamma_k, 1)}\right]+\chi_k\right)M,\quad 1\le k\le q,\lb{3.30}\eea
where $\chi_k=0$ or $1$ for $1\le k\le q$ and $\frac{M\theta}{\pi}\in\Z$
whenever $e^{\sqrt{-1}\theta}\in\sigma(M_k)$ and $\frac{\theta}{\pi}\in\Q$
for some $1\le k\le q$. Furthermore, given $M_0\in\N$, by the proof of Theorem 4.1 of \cite{LoZ}
, we may  further require  $M_0|T$ (since the closure of the set $\{\{Tv\}: T\in\N, \;M_0|T\}$ is still a closed additive subgroup of $\bf T^h$ for some $h\in\N$, where we use notations  as (4.21)  in \cite{LoZ}
.  Then we can use the proof of Step 2 in Theorem 4.1 of \cite{LoZ}
 to get $T$).}

In fact,  Let  $\mu_i=\sum_{\theta\in(0, 2\pi)}S_{M_i}^-(e^{\sqrt{-1}\theta})$ for $1\le i\le q$
and $\alpha_{i, j}=\frac{\theta_j}{\pi}$ where $e^{\sqrt{-1}\theta_j}\in\sigma(M_i)$
for $1\le j\le \mu_i$ and $1\le i\le q$.  Let $h=q+\sum_{1\le i\le q}\mu_i$ and
\bea  v=\left(\frac{1}{M\hat i(\gamma_1, 1)},\dots, \frac{1}{M\hat i(\gamma_q, 1)},
\frac{\alpha_{1, 1}}{\hat i(\gamma_1, 1)}, \frac{\alpha_{1, 2}}{\hat i(\gamma_1, 1)},
\dots\frac{\alpha_{1, \mu_1}}{\hat i(\gamma_1, 1)},
\frac{\alpha_{2, 1}}{\hat i(\gamma_2, 1)},\dots,
\frac{\alpha_{q, \mu_q}}{\hat i(\gamma_q, 1)}\right)\in\R^h.
\lb{3.31}\eea
Then  the above theorem is equivalent to find a vertex 
$$\chi=(\chi_1,\ldots,\chi_q,\chi_{1, 1},\chi_{1, 2},\ldots,\chi_{1, \mu_1},
\chi_{2, 1},\ldots,\chi_{q, \mu_q})$$
of the cube $[0,\,1]^h$ and infinitely many integers $T\in\N$ such that 
\bea|\{Tv\}-\chi|<\epsilon\lb{3.32}\eea
 for any given $\epsilon$ small enough (cf. P. 346 and 349 of \cite{LoZ}
).

{\bf Theorem 3.8.} (cf. Theorem 4.2 of \cite{LoZ}
) {\it Let $H$ be the closure of the subset $\{\{mv\}|m\in\N\}$ in $\bf {T}^h=(\R/\Z)^h$ and
$V=T_0\pi^{-1}H$ be the tangent space of $\pi^{-1}H$ at the origin in $\R^h$,
where $\pi: \R^h\rightarrow  \bf{T}^h$ is the projection map.
Define
\bea A(v)=V\setminus\cup_{v_k\in\R\setminus\Q}\{x=(x_1,,\ldots,x_h)\in V | x_k=0\}.
\lb{3.33}\eea
Define $\psi(x)=0$  when $x\ge 0$ and $\psi(x)=1$ when $x<0$.
Then for any $a=(a_1,\ldots, a_h)\in A(V)$, the vector 
\bea \chi=(\psi(a_1),\ldots,\psi(a_h))\lb{3.34}\eea
makes (\ref{3.32}) hold for infinitely many $T\in\N$. }

{\bf Theorem 3.9.}  (cf. Theorem 4.2 of \cite{LoZ}
) {\it We have the following properties 
for $A(v)$:

(i) When $v\in\R^h\setminus\Q^h$,  then $\dim V\ge 1$,
$0\notin A(v)\subset V$, $A(v)=-A(v)$ and $A(v)$ is open in $V$.

(ii) When $\dim V = 1$, then  $A(v) = V \setminus\{0\}$.

(iii) When $\dim V \ge 2$,  $A(v)$  is obtained from $V$ by deleting all the 
coordinate hyperplanes	with	dimension	
strictly smaller than	$\dim V$	from	$V$.}

\setcounter{equation}{0}
\section{Proof of the main theorems}

In this section, we give the proofs of the main theorems.

 In the rest of this paper, we assume the following:

{\bf (F) There are only finitely many prime closed geodesics
$\{c_j\}_{1\le j\le p}$ on $(M,\,F)$. }

Let $\pi: \tilde M\rightarrow M$ be the universal covering space of $M$
and $\tilde F=\pi^\ast F$, then $(\tilde M, \tilde F)$ is a Finsler manifold 
that is locally isometric to $(M, F)$.
Thus the flag curvature $K$ of $(\tilde M, \tilde F)$ satisfies
$\left(\frac{\lambda}{\lambda+1}\right)^2<K\le 1$ by  assumption.
Hence by \cite{Rad4}, $\tilde M$ is homeomorphic  to $S^n$.

For each  prime closed geodesic $c_j$ on $(M, F)$, clearly there exists a minimal
$\alpha_j\in\N$ such that  $c_j^{\alpha_j}$ lifts to a prime closed geodesic  on 
$(\tilde M, \tilde F)$.  Denote by  $\{\tilde c_{j, 1},\ldots, \tilde c_{j, n_j}\}$
the  lifts of  $c_j^{\alpha_j}$ for some $n_j\in\N$ such that $\tilde c_{j_l}$s are pairwise distinct for $1\le l\le n_j$. For each $l\in \{2,\ldots, n_j\}$, there is a
covering transformation $h:\tilde M\rightarrow \tilde M$ 
such that $h(\tilde c_{j, l})=\tilde c_{j, 1}$.
By the definition of $\tilde F$, $h$ is an isometry on $(\tilde M, \tilde F)$.
Therefore $h$ preserves the energy functional, i.e., $E(\gamma)=E(h(\gamma))$
for any $\gamma\in\Lambda\tilde M$.  In particular,  we have 
\bea  && i(\tilde c_{j, l}^m)=i(\tilde c_{j, 1}^m),\quad \nu(\tilde c_{j, l}^m)=\nu(\tilde c_{j, 1}^m),\nn\\
&&\ol C_q(E, \tilde c_{j,  l}^m) \cong \ol C_q(E,  \tilde c_{j, 1}^m),\quad\forall m\in\N, \;q\in\Z,\;2\le l\le n_j.
\lb{4.1}\eea
Hence there are exactly $q=\sum_{1\le j\le p}n_j$ prime closed geodesics  
$\{\tilde c_{1, 1},\ldots, \tilde c_{1, n_1},\ldots, \tilde c_{p, 1},\ldots, \tilde c_{p, n_p}\}$
on $(\tilde M, \tilde F)$.  Denote by
$\{P_{\tilde c_{j, l}}\}_{1\le j\le p,\,1\le l\le n_j}$ the linearized Poincar\'e maps of 
$\{\tilde c_{j, l}\}_{1\le j\le p,\,1\le l\le n_j}$, then  $P_{\tilde c_{j, l}}=P_{\tilde c_{j, 1}}$
for $1\le j\le p$ and $2\le l\le n_j$.

Since the flag
curvature $K$ of $(\tilde M, \tilde F)$ satisfies
$\left(\frac{\lambda}{\lambda+1}\right)^2<K\le 1$ by assumption,
then every non-constant closed geodesic $\tilde c$ on $(\tilde  M, \tilde F)$ must satisfy
\bea i(\tilde c)\ge n-1, \lb{4.2}\eea
by Theorem 3 and Lemma 3 of \cite{Rad4}.

Now it follows from Theorem 3.6 and (\ref{4.2}) that
\bea i(\tilde c_{j, l}^{m+1})-i(\tilde c_{j, l}^m)-\nu(\tilde c_{j, l}^m)
\ge i(\tilde c_{j, l})-\frac{e(P_{\tilde c_{j, l}})}{2}\ge 0,\quad1\le j\le p,\,1\le l\le n_j,\;\forall m\in\N.\lb{4.3}\eea
Here the last inequality holds by (\ref{4.2}) and the fact that $e(P_{\tilde c_{j, l}})\le 2(n-1)$.

Note that we have $\hat i(\tilde c_{j, l})> n-1$ for $1\le j\le p,\,1\le l\le n_j$
under the pinching assumption by Lemma 2 of \cite{Rad5}.
Hence by Theorem 3.7, (\ref{4.1}) and (\ref{4.2}),
there exist infinitely many
$(N, m_{1, 1},\ldots,m_{1, n_1},\ldots,m_{p, 1},\ldots,m_{p, n_p} )
=(N, m_1,\ldots m_1,\ldots, m_p,\ldots, m_p)\in\N^{q+1}$ such that
\bea
i(\tilde c_{j, l}^{2m_j}) &\ge& 2N-\frac{e(P_{\tilde c_{j, l}})}{2}\ge 2N-(n-1), \lb{4.4}\\
i(\tilde c_{j, l}^{2m_j})+\nu(\tilde c_{j, l}^{2m_j}) &\le& 2N+\frac{e(P_{\tilde c_{j, l}})}{2}\le 2N+(n-1), \lb{4.5}\\
i(\tilde c_{j, l}^{2m_j-m})+\nu(\tilde c_{j, l}^{2m_j-m}) &\le& 2N-(i(\tilde c_{j, l})+2S^+_{P_{\tilde c_{j, l}}}(1)-\nu(\tilde c_{j, l})),\quad \forall m\in\N. \lb{4.6}\\
i(\tilde c_{j, l}^{2m_j+m}) &\ge& 2N+i(\tilde c_{j, l}),\quad\forall m\in\N, \lb{4.7}
\eea
for $1\le j\le p$ and $1\le l\le n_j$.
Moreover $\frac{m_j\theta}{\pi}\in\Z$, whenever $e^{\sqrt{-1}\theta}\in\sigma(P_{\tilde c_{j, l}})$
and $\frac{\theta}{\pi}\in\Q$. In fact, the $m>1$ cases in
(\ref{4.6}) and ({\ref{4.7}) follow from (\ref{4.3}), other parts follow from
Theorem 3.7 directly.
More precisely, by Theorem 3.7
\bea m_j=\left(\left[\frac{N}{M\hat i(\tilde c_{j, l})}\right]+\chi_j\right)M,\quad 1\le j\le p,\lb{4.8}\eea
where $\chi_j=0$ or $1$ for $1\le j\le p$ and $M\in\N$ such that $\frac{M\theta}{\pi}\in\Z$,
whenever $e^{\sqrt{-1}\theta}\in\sigma(P_{\tilde c_{j, l}})$ and $\frac{\theta}{\pi}\in\Q$
for some $1\le j\le p$.

By Theorem 3.5,  there exists a path $f_j\in
C([0,1],\Omega^0(P_{\tilde c_{j, l}}))$ such that $f_j(0)=P_{\tilde c_{j, l}}$ and
\be  f_j(1)= N_1(1,1)^{\diamond p_{j, -}}\diamond I_2^{\diamond p_{j, 0}}\diamond N_1(1,-1)^{\diamond  p_{j, +}}\diamond  G_j,
\qquad1\le j\le p\lb{4.9}\ee
for some nonnegative integers $p_{j, -}$, $p_{j, 0}$, $p_{j, +}$, and some symplectic
matrix $G_j$ satisfying $1\not\in \sigma(G_j)$.
By (\ref{4.9}) and Lemma 3.4 we obtain
\be 2S_{P_{\tilde c_{j, l}}}^+(1) - \nu_1( P_{\tilde c_{j, l}}) = p_{j, -} -p_{j, +} \ge -p_{j, +}\ge 1-n,\quad1\le j\le p,\;1\le l\le n_j. \lb{4.10}\ee
Using (\ref{4.2}) and (\ref{4.10}), the estimates (\ref{4.4})-(\ref{4.7}) become
\bea
i(\tilde c_{j, l}^{2m_j}) &\ge&  2N-(n-1), \lb{4.11}\\
i(\tilde c_{j, l}^{2m_j})+\nu(\tilde c_{j, l}^{2m_j}) &\le& 2N+(n-1), \lb{4.12}\\
i(\tilde c_{j, l}^{2m_j-m})+\nu(\tilde c_{j, l}^{2m_j-m}) &\le& 2N,\quad \forall m\in\N. \lb{4.13}\\
i(\tilde c_{j, l}^{2m_j+m}) &\ge& 2N+(n-1),\quad\forall m\in\N\lb{4.14}
\eea
for $1\le j\le p$ and $1\le l\le n_j$.

In order to prove Theorem 1.2, we need the following:

{\bf Lemma 4.1.} {\it There  exists  $j_0\in\{1,\dots, p\}$ such that
 $i(\tilde c_{j_0, l}^{2m_{j_0}})+\nu(\tilde c_{j_0, l}^{2m_{j_0}})=2N+(n-1)$. Moreover, 
 $\ol{C}_{2N+n-1}(E, \tilde c_{j_0, l}^{2m_{j_0}})\neq 0$ for $1\le l\le n_{j_0}$. }

{\bf Proof.}
We prove  by contradiction, i.e., suppose that 
\bea \ol{C}_{2N+n-1}(E, \tilde c_{j, l}^{2m_{j}})=0,\quad 1\le j\le p,\;1\le l\le n_j.\lb{4.15}\eea
Now  (\ref{4.5})-(\ref{4.7}) becomes 
\bea
i(\tilde c_{j, l}^{m})+\nu(\tilde c_{j, l}^{m}) &\le& 2N\quad\forall m<2m_j, \lb{4.16}\\
i(\tilde c_{j, l}^{2m_j})+\nu(\tilde c_{j, l}^{2m_j}) &\le& 2N+n-1, \lb{4.17}\\
i(\tilde c_{j, l}^{m}) &\ge& 2N+n-1,\quad\forall m>2m_j. \lb{4.18}\eea

By Theorem 3.7,  we can choose $N$ and $\{\chi_j\}_{1\le j\le p}$ such that
\bea \left|\frac{N}{M\hat i(\tilde c_{j, l})}-\left[\frac{N}{M\hat i(\tilde c_{j, l})}\right]-\chi_j\right|<\epsilon
<\frac{1}{1+\sum_{1\le j\le p,\,1\le l\le n_j}4M|\hat\chi(\tilde c_{j, l})|},\quad 1\le j\le p.\lb{4.19}\eea
By Theorem 2.5  we have
\be  \sum_{1\le j\le p, 1\le l\le n_j}\frac{\hat\chi(\tilde c_{j, l})}{\hat{i}(\tilde c_{j, l})}=B(n,1)\in\Q.\lb{4.20}\ee
Note by Theorem 3.7, we can require that $N\in\N$  further satisfies 
\bea 2NB(n,1)\in\Z.\lb{4.21}\eea
Multiplying both sides of (\ref{4.20}) by $2N$ yields
\be  \sum_{1\le j\le p, 1\le l\le n_j}\frac{2N\hat\chi(\tilde c_{j, l})}{\hat{i}(\tilde c_{j, l})}=2NB(n,1).\lb{4.22}\ee

{\bf Claim 1.} {\it We have}
\be  \sum_{1\le j\le p, 1\le l\le n_j}2m_j\hat\chi(\tilde c_{j, l})=2NB(n,1).\lb{4.23}\ee
In fact, by (\ref{4.22}), we have
\bea &&2NB(n,1)\nn\\
=&&\sum_{1\le j\le p, 1\le l\le n_j}\frac{2N\hat\chi(\tilde c_{j, l})}{\hat{i}(\tilde c_{j, l})}\nn\\
=&&\sum_{1\le j\le p, 1\le l\le n_j}2\hat\chi(\tilde c_{j, l})\left(\left[\frac{N}{M\hat i(\tilde c_{j, l})}\right]+\chi_j\right)M
\nn\\&&+\sum_{1\le j\le p, 1\le l\le n_j}2\hat\chi(\tilde c_{j, l})\left(\frac{N}{M\hat i(\tilde c_{j, l})}-\left[\frac{N}{M\hat i(\tilde c_{j, l})}\right]-\chi_j\right)M\nn\\
\equiv&& \sum_{1\le j\le p, 1\le l\le n_j}2m_j\hat\chi(\tilde c_{j, l})+\sum_{1\le j\le p, 1\le l\le n_j}2M\hat\chi(\tilde c_{j, l})\epsilon_j.\lb{4.24}
\eea
By Proposition 2.3 and our choice of $M$, we have
\bea \frac{2m_j}{T(\tilde c_{j, l})}\in\N,\quad 1\le j\le p,\;1\le l\le n_j.\lb{4.25}\eea
Hence (\ref{2.11}) implies that
\bea 2m_j\hat\chi(\tilde c_{j, l})\in\Z,\quad1\le j\le p,\,1\le l\le n_j.\lb{4.26}\eea
Now Claim 1 follows by (\ref{4.19}), (\ref{4.21}), (\ref{4.24}) and  (\ref{4.26}).

{\bf Claim 2.} {\it  We have
\be \sum_{1\le j\le p, 1\le l\le n_j}2m_j\hat\chi(\tilde c_{j, l})=M_0-M_1+M_2-\cdots+(-1)^{2N+n-2}M_{2N+n-2}.\lb{4.27}\ee}
In fact, by definition, the right hand side of (\ref{4.27}) is
\bea RHS=\sum_{q\le 2N+n-2,\, m\ge 1\atop 1\le j\le p,\,1\le l\le n_j}(-1)^q\dim\ol{C}_q(E, \tilde c_{j, l}^m).\lb{4.28}
\eea
By (\ref{4.15})-(\ref{4.18}) and Proposition 2.1, we have
\bea RHS&=&\sum_{ q\le 2N+n-2,\,1\le m\le 2m_j\atop1\le j\le p,\,1\le l\le n_j}(-1)^q\dim\ol{C}_q(E, \tilde c_{j, l}^m),\lb{4.29}\\
&=&\sum_{1\le m\le 2m_j\atop 1\le j\le p,\;1\le l\le n_j}\chi(\tilde c_{j, l}^m),\lb{4.30}
\eea
where the second equality follows from (\ref{2.9}), (\ref{4.15})-(\ref{4.17})
and Proposition 2.1.

By Proposition 2.3, (\ref{2.9})-(\ref{2.11}) and (\ref{4.25}), we have
\bea  \sum_{1\le m\le 2m_j }\chi(\tilde c_{j, l}^m)
&=&\sum_{0\le s< 2m_j/T(\tilde c_{j, l})\atop1\le m\le T(\tilde c_{j, l})}\chi(\tilde c_{j, l}
^{sT(\tilde c_{j, l}
)+m})\nn\\
&=&\frac{2m_j}{T(\tilde c_{j, l}
)}\sum_{1\le m\le T(\tilde c_{j, l}
)}\chi(\tilde c_{j, l}
^m)\nn\\
&=& 2m_j\hat\chi(\tilde c_{j, l}
),\lb{4.31}\eea
This proves Claim 2.

In order to prove the lemma, we  consider the following two
cases according to the parity of $n$.

{\bf Case 1. } $n=2k+1$ is odd.

In this case, we have by (\ref{2.13})
\be B(n,1)=\frac{n+1}{2(n-1)}=\frac{k+1}{2k}.\lb{4.32}\ee
By Theorem 3.7 we may  further assume  $N=mk$ for some $m\in\N$.

Thus by (\ref{4.23}), (\ref{4.27}) and (\ref{4.32}), we have
\be M_0-M_1+M_2-\cdots+(-1)^{2N+n-2}M_{2N+n-2}=m(k+1).\lb{4.33}
\ee
On the other hand, we have by (\ref{2.15})
\bea &&b_0-b_1+b_2-\cdots+(-1)^{2N+n-2}b_{2N+n-2}\nn\\
=&&b_{2k}+(b_{2k+2}+\cdots+b_{4k}+\cdots+b_{2mk+2}+\cdots+b_{2mk+2k})-b_{2mk+2k}\nn\\
=&& 1+m(k-1+2)-2\nn\\
=&& m(k+1)-1.\lb{4.34}
\eea
In fact, we cut off the sequence $\{b_{2k+2},\ldots,b_{2mk+2k}\}$ into
$m$ pieces, each of them contains  $k$ terms. Moreover, each piece contain
$1$ for $k-1$ times and  $2$ for one time. Thus (\ref{4.34}) holds.

Now by Theorem 2.7 and (\ref{4.34}), we have
\bea -m(k+1)&=&M_{2N+n-2}-M_{2N+n-3}+\cdots+M_1-M_0\nn\\
&\ge&b_{2N+n-2}-b_{2N+n-3}+\cdots+b_1-b_0\nn\\
&=&-(m(k+1)-1).\lb{4.35}
\eea
This contradiction yields the lemma for $n$ being odd.

{\bf Case 2. } $n=2k$ is even.

In this case, we have by (\ref{2.13})
\be B(n,1)=\frac{-n}{2n-2}=\frac{-k}{2k-1}.\lb{4.36}\ee
As in Case 1, we may assume $N=m(2k-1)$ for some $m\in\N$.

Thus by (\ref{4.23}), (\ref{4.27}) and (\ref{4.36}), we have
\be M_0-M_1+M_2-\cdots+(-1)^{2N+n-2}M_{2N+n-2}=-2mk.\lb{4.37}
\ee
On the other hand, we have by (\ref{2.17})
\bea &&b_0-b_1+b_2-\cdots+(-1)^{2N+n-2}b_{2N+n-2}\nn\\
=&&-b_{2k-1}-(b_{2k+1}+\cdots+b_{6k-3}+\cdots+b_{(m-1)(4k-2)+2k+1}+\cdots+b_{m(4k-2)+2k-1})\nn\\
&&+b_{m(4k-2)+2k-1}\nn\\
=&& -1-m(2k-2+2)+2\nn\\
=&& -2mk+1.\lb{4.38}
\eea
In fact, we cut off the sequence $\{b_{2k+1},\ldots,b_{m(4k-2)+2k-1}\}$ into
$m$ pieces, each of them contains  $2k-1$ terms. Moreover, each piece contain
$1$ for $2k-2$ times and  $2$ for one time. Thus (\ref{4.38}) holds.

Now by (\ref{4.37}), (\ref{4.38}) and Theorem 2.7, we have
\bea -2mk&=&M_{2N+n-2}-M_{2N+n-3}+\cdots+M_1-M_0\nn\\
&\ge&b_{2N+n-2}-b_{2N+n-3}+\cdots+b_1-b_0\nn\\
&=&-2mk+1.\lb{4.39}
\eea
This contradiction yields the lemma for $n$ being even.
 \hfill\hb

We will use the following theorem of N. Hingston. Note that the proof of N. Hingston's theorem does
not use the special properties of Riemannian metric, hence it
holds for Finsler metric as well.

{\bf Theorem 4.2.} (Follows from Proposition 1 of \cite{Hin2}, cf.
Lemma 3.4.12 of \cite{Kli3})
{\it Let $c$ be a closed geodesic of length $L$ on a compact Finsler
manifold $(M, F)$ such that as a critical orbit of the
energy functional $E$ on $\Lm M$, every orbit $S^1\cdot c^m$ of
its iteration $c^m$ is isolated. Suppose
\bea
&& i(c^m)+\nu(c^m) \le m(i(c)+\nu(c))-(n-1)(m-1),  \qquad
                \forall m\in\N, \lb{4.40}\\
&& k_{\nu(c)}(c)\neq0.  \lb{4.41} \eea
Then  $(M, F)$ has infinitely many prime closed geodesics. }

Note that in (\ref{4.41}), we have used the Shifting theorem in
\cite{GrM1}. Especially, (\ref{4.41}) means that $c$ is a local
maximum in the local characteristic manifold $N_c$ at $c$.

{\bf Lemma 4.3.} {\it The closed geodesic $\tilde c_{j_0, l}$ found in
Lemma 4.1 satisfy the following: 

(i) $e(P_{\tilde c_{j_0, l}})=2n-2$, i.e., $\tilde c_{j_0, l}$ is elliptic.

(ii) $P_{\tilde c_{j_0, l}}$ does not contain $N_1(1,1)$, $N_1(-1,-1)$  and nontrivial $N_2(\omega, b)$.

(iii) Any trivial $N_2(\omega, b)$ contained in  $P_{\tilde c_{j_0, l}}$ must satisfies
$\frac{\theta}{\pi}\in\Q$, where $\omega=e^{\sqrt{-1}}\theta$.

(iv) $k_{\nu({\tilde c_{j_0, l}}^{T({\tilde c_{j_0, l}})})}({\tilde c_{j_0, l}}^{T({\tilde c_{j_0, l}})})\neq 0$.
Hence ${\tilde c_{j_0, l}}^{T({\tilde c_{j_0, l}})}$ is a local maximum of the energy functional
in the local characteristic manifold at $\tilde c_{j_0, l}^{T({\tilde c_{j_0, l}})}$.

(v) $P_{\tilde c_{j_0, l}}$ must contain a term $R(\theta)$ with $\frac{\theta}{\pi}\notin\Q$.
}

{\bf Proof.} Note that by (\ref{3.26}), Lemma 4.1 and Theorem 3.6, we have 
\bea 2N+(n-1)&=&i(\tilde c_{j_0, l}^{2m_{j_0}})+\nu(\tilde c_{j_0, l}^{2m_{j_0}})
\nn\\&\le& i(\tilde c_{j_0, l}^{2m_{j_0}+1})-i(\tilde c_{j_0, l})+\frac{e(P_{\tilde c_{j_0, l}})}{2}\lb{4.42}\\
&=& 2N+\frac{e(P_{\tilde c_{j_0, l}})}{2}\le 2N+(n-1).\lb{4.43}\eea
Hence (i) holds. If any one of (ii)-(iii) does not hold, then by Theorem 3.5, the strict inequality in
(\ref{4.42}) must hold and yields a contraction. 

(iv) follows directly from Proposition 2.3 and Lemma 4.1.

We prove (v) by contradiction. 
Consider $g_l=\tilde c_{j_0, l}^{T(\tilde c_{j_0, l})}$. Then by (i)-(iii) and the assumption,
$P_{g_l}$ can be connected in
$\Omega^0(P_{g_l})$ to $I_{2p_0^\prime}\diamond N_1(1,\,-1)^{\diamond p_+^\prime}$
with $p_0^\prime+p_+^\prime=n-1$ as in Theorem 3.5.
In fact, by (i)-(iii)  and the assumption, the basic normal
form decomposition (\ref{3.19}) in Theorem 3.5 becomes
\bea P_{\tilde c_{j_0, l}}=&&I_{2p_0}\diamond N_1(1,\,-1)^{\diamond p_+}
\diamond N_1(-1,\,1)^{\diamond q_-}\diamond(-I_{2q_0})\nn\\
&&\diamond R(\theta_1)\diamond\cdots\diamond R(\theta_r)
\diamond N_2(\lambda_1,\,v_1)\diamond\cdots\diamond N_2(\lambda_{r_0},\,v_{r_0})
\nn\eea
together with $\frac{\theta_j}{\pi}\in\Q$ for
$1\le j\le r$, $\frac{\beta_j}{\pi}\in\Q$  for $1\le j\le r_0$
and $p_0+p_++q_-+q_0+r+2r_0=n-1$.
By Proposition 2.3, we have $\frac{T(\tilde c_{j_0, l})\theta_j}{2\pi}\in\Z$,
$\frac{T(\tilde c_{j_0, l})\beta_j}{2\pi}\in\Z$  and $2|T(\tilde c_{j_0, l})$ whenever
$-1\in\sigma(P_{\tilde c_{j_0, l}})$. Hence $R(\theta_j)^{T(\tilde c_{j_0, l})}=I_2$,
$N_2(\lambda_j,\,v_j)^{T(\tilde c_{j_0, l})}$ can be connected within
$\Omega^0(N_2(\lambda_j,\,v_j)^{T(\tilde c_{j_0, l})})$ to $N_1(1,\,-1)^{\diamond 2}$
and $(-I_2)^{T(\tilde c_{j_0, l})}=I_2$, $N_1(-1, \,1)^{T(\tilde c_{j_0, l})}$
can be connected within $\Omega^0(N_1(-1, \,1)^{T(\tilde c_{j_0, l})})$ to $N_1(1,\,-1)$
whenever $-1\in\sigma(P_{\tilde c_{j_0, l}})$. Thus $p_0^\prime=p_0+q_0+r$ and
$p_+^\prime=p_++q_-+2r_0$ and then $P_{g_l}$ behaves as claimed.

Now by Theorem 3.5, we have
\bea i(g_l^m)=m(i(g_l)+p_0^\prime)-p_0^\prime,\quad \nu(g_l^m)\equiv2p_0^\prime+p_+^\prime.\qquad\forall m\in\N.\lb{4.44}
\eea
Hence
\bea i(g_l^m)+\nu(g_l^m)=m(i(g_l)+p_0^\prime)+p_0^\prime+p_+^\prime.\qquad\forall m\in\N.\lb{4.45}
\eea
On the other hand
\bea &&m(i(g_l)+\nu(g_l))-(n-1)(m-1)\nn\\
=&&m(i(g_l)+2p_0^\prime+p_+^\prime)-(p_0^\prime+p_+^\prime)(m-1)\nn\\
=&&m(i(g_l)+p_0^\prime)+p_0^\prime+p_+^\prime.\qquad\forall m\in\N.\lb{4.46}
\eea
By  (iv), we have
\be k_{\nu(g_l)}(g_l)=k_{\nu({\tilde c_{j_0, l}}^{T({\tilde c_{j_0, l}})})}({\tilde c_{j_0, l}}^{T({\tilde c_{j_0, l}})})\neq 0.\lb{4.47}
\ee
Hence we can use Theorem 4.2 to obtain infinitely many prime closed geodesics,
which contradicts to the assumption (F). This complete the proof of Lemma 4.3.\hfill\hb

{\bf Lemma 4.4.} {\it There exists no closed geodesic $g$ on $(\tilde M, \tilde F)$
such that  $P_g=N_1(1, -1)^{\diamond (n-1)}$ and $k_{\nu(g)}(g)\neq 0$.}

{\bf Proof.} By the proof of Theorem 4.3,  we can use Theorem 4.2 to obtain infinitely many prime closed geodesics,
which contradicts to the assumption (F). \hfill\hb

{ \bf Proof of Theorem 1.2.}  
By Lemma 2.10, for every $i\in\N$, there exist some $m(i), j(i)\in\N$
such that
\bea E(\tilde c_{j(i), l}^{m(i)})=\kappa_i,\quad
\ol{C}_{2i+\dim(z)-2}(E, \tilde c_{j(i), l}^{m(i)})\neq 0,\lb{4.48}\eea
and by \S2, we have $\dim(z)=n+1$.

{\bf Claim 1.} {\it We have the following}
\bea m(i)=2m_{j(i)},  \quad{\rm if}\quad  2i+\dim(z)-2\in (2N,\,2N+n-1),
\lb{4.49}\eea
In fact, we have
\bea \ol{C}_q(E, \tilde c_{j, l}^m)= 0,\qquad {\rm if}\quad q\in (2N,\,2N+n-1)
\lb{4.50}\eea
for $1\le j\le p$, $1\le l\le n_j$ and $m\neq 2m_j$ by (\ref{4.13}), (\ref{4.14}) and  Proposition 2.1.
Thus in order to satisfy (\ref{4.48}), we must have $m(i)=2m_{j(i)}$.

By Lemma 2.9,  for $2\alpha+\dim(z)-2,  2\beta+\dim(z)-2 \in (2N,\,2N+n-1)$ with 
 $\alpha\neq \beta$
\bea  E(\tilde c_{j(\alpha), l}^{2m_{j(\alpha)}})=\kappa_\alpha
\neq \kappa_\beta= E(\tilde c_{j(\beta), k}^{2m_{j(\beta)}}).
\lb{4.51}\eea
Thus $\tilde c_{j(\alpha), l}\neq \tilde c_{j(\beta), k}$ for $1\le l\le n_{j(\alpha)}$
and $1\le k\le n_{j(\beta)}$.
Therefore there are
\be ^\#\{i: 2i+\dim(z)-2\in (2N,\,2N+n-1)\}=\left[\frac{n}{2}\right]-1
\lb{4.52}\ee
closed geodesics on $(M, F)$. By a permutation of $\{1,\ldots, p\}$,
we may denote these closed geodesics by $\{c_1,\ldots, c_{\left[\frac{n}{2}\right]-1}\}$.

{\bf Claim 2.} {\it  If $n$ is odd, then it is impossible that 
$\ol{C}_{2N}(E, \tilde c_{j, l}^{2m_j-m})\neq 0$ for some $1\le j\le p$ and $m\in\N$.}

Suppose the contrary, then by (\ref{4.3}), (\ref{4.6}), (\ref{4.10})  and Proposition 2.1, 
we have  $i(\tilde c_{j, l}^{2m_j-1})+\nu(\tilde c_{j, l}^{2m_j-1})=2N$,
$P_{\tilde c_{j, l}}=N_1(1, -1)^{\diamond (n-1)}$ together with 
$i(\tilde c_{j, l})=n-1$ for some $1\le j\le p$.
Thus by Proposition 2.3,  we have $k_{\nu(\tilde c_{j, l})}(\tilde c_{j, l})\neq 0$.
Hence we can use Lemma 4.4 to obtain infinitely many prime closed geodesics,
which contradicts to the assumption (F).  This proves Claim 2.

Thus by Claim 2, (\ref{4.14}) and Proposition 2.1,  for $n$ being odd and $2i+\dim (z)-2=2N$
we have 
\bea E(\tilde c_{j(i), l}^{2m_{j(i)}})=\kappa_i,\qquad \ol{C}_{2N}(E, \tilde c_{j(i), l}^{2m_{j(i)}})\neq 0,
\lb{4.53}\eea
 thus we have 
one more closed geodesic on $(M, F)$ by Lemma 2.9. 
Hence we have $[\frac{n+1}{2}]-1$ closed geodesics on $(M, F)$.
We may denote these closed geodesics by $\{c_1,\ldots, c_{\left[\frac{n+1}{2}\right]-1}\}$.

By Lemma 4.1, there  exists  $j_0\in\{1,\dots, p\}$ such that
 $i(\tilde c_{j_0, l}^{2m_{j_0}})+\nu(\tilde c_{j_0, l}^{2m_{j_0}})=2N+(n-1)$, moreover, 
 $\ol{C}_{2N+n-1}(E, \tilde c_{j_0, l}^{2m_{j_0}})\neq 0$ for $1\le l\le n_{j_0}$. 
 By Proposition 2.2, we have
 \bea \ol{C}_{q}(E, \tilde c_{j_0, l}^{2m_{j_0}})=0,\quad q\neq 2N+(n-1),\;1\le l\le n_{j_0}.  \lb{4.54}\eea
 By (\ref{4.48}), (\ref{4.53}) and (\ref{4.54}),  $c_{j_0}\notin\{c_1,\ldots, c_{\left[\frac{n+1}{2}\right]-1}\}$.
Hence we have $[\frac{n+1}{2}]$ closed geodesics on $(M, F)$.
The proof of  Theorem 1.2  is complete.\hfill\hb

{\bf Proof of Theorem 1.3.} By Theorem 1.1 of \cite{Liu1} or Lemma 3.2 of \cite{LiX}, $c_j$ is hyperbolic if and only if 
$\tilde c_{j ,l}$ are hyperbolic for $1\le l\le n_j$.
By  (\ref{4.4}) and (\ref{4.5}), a hyperbolic closed geodesic 
$\tilde c_{j, l}$ must satisfy $i(\tilde c_{j, l}^{2m_j})=2N=i(\tilde c_{j, l}^{2m_j})+\nu(\tilde c_{j, l}^{2m_j})$. Thus by Proposition 2.1 we have 
 \bea \ol{C}_{q}(E, \tilde c_{j, l}^{2m_{j}})=0,\quad q\neq 2N.\lb{4.55}\eea
By (\ref{4.48}) and (\ref{4.55}), the closed geodesics  $\{c_1,\ldots, c_{\left[\frac{n}{2}\right]-1}\}$ are non-hyperbolic.
By Lemma 4.1 and (\ref{4.55}), the closed geodesic $c_{j_0}$ is non-hyperbolic.
Therefore there are $[\frac{n}{2}]$ non-hyperbolic closed geodesics  on $(M, F)$.\hfill\hb

{\bf Proof of Theorem 1.4.}  By Theorem 1.1 of \cite{Liu1} or Lemma 3.2 of \cite{LiX}, $P_{c_j}$ possess an eigenvalue of the form $\exp(\pi i \mu)$ with an
irrational $\mu$
 if and only if 
$P_{\tilde c_{j ,l}}$ have the same property for $1\le l\le n_j$.

Firstly we prove the following:

{\bf Claim 1.} {\it There are at least $[\frac{n-1}{2}]$ closed geodesics 
$c_{j_k}$ for $1\le k\le [\frac{n-1}{2}]$ on $(M, F)$ 
 such that $\tilde c_{j_k, l}\in \mathcal{V}_\infty(\tilde M, \tilde F)$ for $1\le k\le [\frac{n-1}{2}]$ and $1\le l\le n_{j_k}$. }

As in the proof of Theorem 1.2, for any $N$ chosen in (\ref{4.11})-(\ref{4.14}) fixed
and $2i+\dim(z)-2\in [2N,\,2N+n-1)$, there exist some $1\le j(i)\le p$
such that $\tilde c_{j(i)}$ is $(2m_{j(i)}, i)$-variationaly visible
by (\ref{4.48}), (\ref{4.49})  and  (\ref{4.53}).
Moreover, if $i_1\neq i_2$, then we must have $j(i_1)\neq j(i_2)$
by (\ref{4.51}), (\ref{4.53}) and Lemma 2.9.
Hence the map
\be \Psi: (2\N+\dim(z)-2)\cap [2N,\,2N+n-1)\rightarrow\{c_j\}_{1\le j\le p},
\qquad 2i+\dim(z)-2\mapsto c_{j(i)}
\lb{4.56}\ee
is injective. We remark here that if there are more that one $c_j$
satisfy (\ref{4.48}), we take any one of it. 
Since we have infinitely many $N$ satisfying (\ref{4.11})-(\ref{4.14})
and  the number of prime closed geodesics is finite, Claim 1 must hold.

{\bf Claim 2.} {\it Among the $[\frac{n-1}{2}]$ closed geodesics $\{c_{j_k}\}_{1\le k\le [\frac{n-1}{2}]}$  in Claim 1, there are at least $[\frac{n-1}{2}]-1$ ones
possessing irrational average indices.}

We prove the claim as the following: By Theorem 3.7, we can obtain
infinitely many $N$ in (\ref{4.11})-(\ref{4.14}) satisfying the
further properties:
\be\frac{N}{M\hat i(\tilde c_{j, l})}\in\N\quad{\rm and}\quad\chi_j=0, \qquad{\rm if}\quad \hat i(\tilde c_{j, l})\in\Q.\lb{4.57}
\ee
Now suppose $\hat i(\tilde c_{j_{k_1}, l})\in\Q$ and
$\hat i(\tilde c_{j_{k_2}, l})\in\Q$ hold for some distinct $1\le k_1, k_2\le  [\frac{n-1}{2}]-1$.
Then by (\ref{4.8}) and (\ref{4.57}) we have
\bea  2m_{j_{k_1}}\hat i(\tilde c_{j_{k_1}, l})&=&2\left(\left[\frac{N}{M\hat i(\tilde c_{j_{k_1}, l})}\right]+\chi_{j_{k_1}}\right)M\hat i(\tilde c_{j_{k_1}, l})
\nn\\&=&2\left(\frac{N}{M\hat i(\tilde c_{j_{k_1}, l})}\right)M\hat i(\tilde c_{j_{k_1}, l})=2N
=2\left(\frac{N}{M\hat i(\tilde c_{j_{k_2}, l})}\right)M\hat i(\tilde c_{j_{k_2}, l})
\nn\\&=&2\left(\left[\frac{N}{M\hat i(\tilde c_{j_{k_2}, l})}\right]+\chi_{j_{k_2}}\right)M\hat i(\tilde c_{j_{k_2}, l})
=2m_{j_{k_2}}\hat i(\tilde c_{j_{k_2}, l}).
\lb{4.58}\eea
On the other hand, by (\ref{4.56}), we have
\be \Psi(2i_1+\dim(z)-2)=c_{j_{k_1}},\quad\Psi(2i_2+\dim(z)-2)=c_{j_{k_2}},\qquad
{\rm for \;some}\quad i_1\neq i_2.\lb{4.59}\ee
Thus by (\ref{4.48}), (\ref{4.49}),  (\ref{4.53}) and Lemma 2.9, we have
\be E(\tilde c_{j_{k_1}, l}^{2m_{j_{k_1}}})=\kappa_{i_1}\neq\kappa_{i_2}=E(\tilde c_{j_{k_2}, l}^{2m_{j_{k_2}}}).\lb{4.60}\ee
Since $\tilde c_{j_{k_1}, l}, \tilde c_{j_{k_2}, l}\in\mathcal{V}_\infty(\tilde M, \tilde F)$, by Theorem 2.12 we have
\be \frac{\hat i(\tilde c_{j_{k_1}, l})}{L(\tilde c_{j_{k_1}, l})}=2\sigma=\frac{\hat i(\tilde c_{j_{k_2}, l})}{L(\tilde c_{j_{k_2}, l})}.\lb{4.61}
\ee
Note that we have the relations
\be L(\tilde c^m)=mL(\tilde c),\quad \hat i(\tilde c^m)=m\hat i(\tilde c),\quad L(\tilde c)=\sqrt{2E(\tilde c)},\qquad\forall m\in\N,\lb{4.62}\ee
for any closed geodesic $\tilde c$ on $(\tilde M,\,\tilde  F)$.

Hence we have
\bea 2m_{j_{k_1}}\hat i(\tilde c_{j_{k_1}, l})
&=&2\sigma\cdot 2m_{j_{k_1}}L(\tilde c_{j_{k_1}, l})
=2\sigma L(\tilde c_{j_{k_1}, l}^{2m_{j_{k_1}}})\nn\\
&=&2\sigma\sqrt{2E(\tilde c_{j_{k_1}, l}^{2m_{j_{k_1}}})}
=2\sigma\sqrt{2\kappa_{i_1}}\nn\\
&\neq&2\sigma\sqrt{2\kappa_{i_2}}
=2\sigma\sqrt{2E(\tilde c_{j_{k_2}, l}^{2m_{j_{k_2}}})}\nn\\
&=&2\sigma L(\tilde c_{j_{k_2}, l}^{2m_{j_{k_2}}})
=2\sigma\cdot 2m_{j_{k_2}}L(\tilde c_{j_{k_2}, l})
=2m_{j_{k_2}}\hat i(\tilde c_{j_{k_2}, l}).
\lb{4.63}\eea
This contradict to (\ref{4.58}) and then we must have
$\hat i(\tilde c_{j_{k_1}, l})\in\R\setminus\Q$ or $\hat i(\tilde c_{j_{k_2}, l})\in\R\setminus\Q$.
Hence there is at most one $1\le k\le  [\frac{n-1}{2}]$
such that $\hat i(\tilde c_{j_{k}, l})\in\Q$, i.e., there are at least $[\frac{n-1}{2}]-1$ ones
possessing irrational average indices.
This proves Claim 2.

Suppose $c_{j_k}$ is any closed geodesic found in Claim 2, 
then we have $\hat i(\tilde c_{j_k, l})\in\R\setminus\Q$.
Therefore by Theorem 3.5, the linearized Poincar\'e map $P_{\tilde c_{j_k, l}}$
of $\tilde c_{j_k, l}$ must contains a  term $R(\theta)$ with $\frac{\theta}{\pi}\notin\Q$.
By  Lemma 4.3, the closed geodesic  $\tilde c_{j_0, l}$ also has this property.
Moreover, we have $c_{j_0}\notin \{c_{j_k}\}_{1\le k\le [\frac{n-1}{2}]}$
as in the proof of Theorem 1.2. Thus Theorem 1.4 is true.\hfill\hb

{\bf Proof of Theorem 1.5.} 
We prove by contradiction. i.e., assume there are exactly two closed geodesics $c_1$ and $c_2$
on $(M, F)$ by Theorem 1.2.

{\bf Step 1.} {\it We can write
\bea P_{\tilde c_{j, l}}=R(\theta_j)\diamond M_j^\prime,\qquad 1\le  j\le 2,\;1\le  l\le n_j,\lb{4.64}
\eea
with  $\frac{\theta_j}{\pi}\notin\Q$ and $M_j^\prime\in \Sp(2)$.
}

By Lemma 4.3,   we may assume 
$P_{\tilde c_{1, l}}$  contains a term $R(\theta_1)$ with $\frac{\theta_1}{\pi}\notin\Q$.
and  
\bea k_{\nu({\tilde c_{1, l}}^{T({\tilde c_{1, l}})})}({\tilde c_{1, l}}^{T({\tilde c_{1, l}})})\neq 0.
\lb{4.65}\eea
Thus in order to prove Step 1, we only need to prove (\ref{4.64}) for $j=2$.

By Theorem 3.5, we have $\hat i(\tilde c_{1, l})=i(\tilde c_{1, 1})+p_-+p_0-r+\sum_{1\le j\le r}\frac{\phi_j}{\pi}$ for some $p_-, p_0, r\in\N_0$.
By Lemma 4.3, we have $r\ge 1$ and then 
\bea \hat i(\tilde c_{1, l})=i(\tilde c_{1, 1})+p_-+p_0-r+\sum_{2\le j\le r}\frac{\phi_j}{\pi}+\frac{\theta_1}{\pi}\equiv\Delta+\frac{\theta_1}{\pi}.\lb{4.66}\eea
and $\frac{\alpha_{1, 1}}{\hat i(\tilde c_{1, l})}=\frac{\theta_1/\pi}{\Delta+\theta_1/\pi}$,
where we use notations as in Theorem 3.7, i.e., $\alpha_{1, 1}=\frac{\theta_1}{\pi}$.
We have the following three cases:

{\bf Case 1.1.} {\it We have $\Delta\in\Q$.}

Denote by $\beta=\Delta+\frac{\theta_1}{\pi}\notin\Q$. Then we have
$(\frac{1}{M\hat i(\tilde c_{1, l})},\, \frac{\alpha_{1, 1}}{\hat i(\tilde c_{1, l})})=(\frac{1}{M\beta}, \, 1-\frac{\Delta}{\beta})$.
Thus if $\frac{N}{M\hat i(\tilde c_{1, l})}=K+\mu$ for some $K\in\Z$ and $\mu\in (-1,\,1)$,
we have $\frac{N\alpha_{1, 1}}{\hat i(\tilde c_{1, l})}=N-M\Delta K-M\Delta\mu$.
Note that by the choice of $M$, we have $M\Delta\in\Z$, 
therefore by (\ref{3.31}) and  (\ref{3.32}), we have
\be
 \left\{\matrix{\chi_{1, 1}=1\quad{\rm if}\quad \chi_1=0, \cr
\chi_{1, 1}=0\quad{\rm if}\quad \chi_1=1. \cr}\right.  \lb{4.67}\ee
Thus either $(\chi_1,\,\chi_{1, 1})=(1,\, 0)$  or $(\chi_1,\,\chi_{1, 1})=(0,\, 1)$
holds. By (4.16)  and (4.17) of \cite{LoZ},
we have
\bea
 &&\{m_1\alpha_{1, 1}\}=\left\{\left\{\frac{N\alpha_{1, 1}}{\hat i(\tilde c_{1, l})}\right\}-\chi_{1, 1}
+\left(\chi_1-\left\{\frac{N}{M\hat i(\tilde c_{1, l})}\right\}\right)M\alpha_{1, 1}\right\}
=\{A_{1, 1}(N)+B_{1, 1}(N)\}
\nn\\=&& \left\{\matrix{\left\{\left\{\frac{N\alpha_{1, 1}}{\hat i(\tilde c_{1, l})}\right\}-\chi_{1, 1}
+\left(\chi_1-\left\{\frac{N}{M\hat i(\tilde c_{1, l})}\right\}\right)M\alpha_{1, 1}\right\}
\quad{\rm if}\quad (\chi_1,\,\chi_{1, 1})=(1,\, 0), \cr
\left\{1+\left\{\frac{N\alpha_{1, 1}}{\hat i(\tilde c_{1, l})}\right\}-\chi_{1, 1}
+\left(\chi_1-\left\{\frac{N}{M\hat i(\tilde c_{1, l})}\right\}\right)M\alpha_{1, 1}\right\}
\quad{\rm if}\quad (\chi_1,\,\chi_{1, 1})=(0,\, 1), \cr}\right.
\nn\eea
where $A_{1, 1}(N)=\left\{\frac{N\alpha_{1, 1}}{\hat i(\tilde c_{1, l})}\right\}-\chi_{1, 1}$
and $B_{1, 1}(N)=\left(\chi_1-\left\{\frac{N}{M\hat i(\tilde c_{1, l})}\right\}\right)M\alpha_{1, 1}$.
In fact, we have  $A_{1, 1}(T)>0$,  $B_{1, 1}(T)>0$ for  $(\chi_1,\,\chi_{1, 1})=(1,\, 0)$,
and    $A_{1, 1}(T)<0$,  $B_{1, 1}(T)<0$ for  $(\chi_1,\,\chi_{1, 1})=(0,\, 1)$,
thus the last equality above holds.

Hence by (\ref{3.32}), we have
\bea\left\{\matrix{\{m_1\alpha_{1, 1}\}<(2M+1)\epsilon\quad&&{\rm if}\quad (\chi_1,\,\chi_{1, 1})=(1,\, 0), \cr
\{m_1\alpha_{1, 1}\}>1-(2M+1)\epsilon\quad&&{\rm if}\quad (\chi_1,\,\chi_{1, 1})=(0,\, 1), \cr}\right.
\lb{4.68}\eea
where we have used the fact that $\alpha_{1, 1}=\theta_1/\pi\in(0,\,2)$.

By choosing $\epsilon\in \left(0,\;\frac{1}{2M+1}\min\{\frac{\theta_1}{2\pi},\,1-\frac{\theta_1}{2\pi}\}\right)$, by Theorem 3.5, we have
\bea &&i(\tilde c_{1, l}^{2m_1+1})-i(\tilde c_{1, l}^{2m_1})\nn\\
=&&i(\tilde c_{1, l})+p_-+p_0-r+\sum_{1\le i\le r}\left(2E\left(\frac{(2m_1+1)\phi_i}{2\pi}\right)-2E\left(\frac{2m_1\phi_i}{2\pi}\right)\right)+q_0+q_+\nn\\=&&2E\left(\frac{(2m_1+1)\theta_1}{2\pi}\right)-2E\left(\frac{2m_1\theta_1}{2\pi}\right)+i(\tilde c_{1, l})+p_-+p_0-r+2(r-1)+q_0+q_+\nn\\
\equiv&&2E\left(\frac{(2m_1+1)\theta_1}{2\pi}\right)-2E\left(\frac{2m_1\theta_1}{2\pi}\right)+\Pi\nn\\=&&2\left(E\left(\frac{2m_1\theta_1}{2\pi}+\frac{\theta_1}{2\pi}\right)
-E\left(\frac{2m_1\theta_1}{2\pi}\right)\right)+\Pi
\nn\\=&&2\left(E\left(\{m_1\alpha_{1, 1}\}+\frac{\theta_1}{2\pi}\right)
-E\left(\{m_1\alpha_{1, 1}\}\right)\right)+\Pi
\nn\\=&&\left\{\matrix{\Pi\quad{\rm if}\quad \{m_1\alpha_{1, 1}\}<(2M+1)\epsilon,\qquad \cr
2+\Pi\quad{\rm if}\quad \{m_1\alpha_{1, 1}\}>1-(2M+1)\epsilon, \cr}\right.
\lb{4.69}\eea
where   $p_-, p_0, q_0, q_+\in \N_0$,  $r\in\N$, $\phi_1=\theta_1$ and 
$\frac{\phi_2}{2\pi}\in(0, 1)\cap\Q$, $\frac{m_1\phi_2}{\pi}\in\Z$ whenever $r=2$.

Hence by Theorems 3.7-3.9, we can choose another $N^\prime\in\N$ 
and
\bea m_j^\prime=\left(\left[\frac{N^\prime}{M\hat i(\tilde c_{j, l})}\right]+\chi_j^\prime\right)M,\quad 1\le j\le 2,\lb{4.70}\eea
such that   $\{m_1\alpha_{1, 1}\}<(2M+1)\epsilon$
and $ \{m_1^\prime\alpha_{1, 1}\}>1-(2M+1)\epsilon$.
By (\ref{3.26}), (\ref{4.69})  and $\nu(\tilde c_{1, l}^{2m_1})=\nu(\tilde c_{1, l}^{2m_1^\prime})$, we have 
\bea i(\tilde c_{1, l}^{2m_1^\prime})+ \nu(\tilde c_{1, l}^{2m_1^\prime})
\le 2N^\prime+(n-1)-2.\lb{4.71}\eea
Hence  the closed geodesic found in Lemma 4.1 for $N^\prime$ must be $\tilde c_{2, l}$ by Proposition 2.1,
and then Step 1 holds in this case by Lemma 4.3.

{\bf Case 1.2.}  {\it We have $\Delta\notin\Q$ and $\hat i(\tilde c_{1, l})\in\Q$.} 

In this case we have $P_{\tilde c_{1, l}}=R(\phi_1)\diamond R(\phi_2)$
with $\frac{\phi_i}{\pi}\notin\Q$ for $i=1, 2$ and $\phi_1=\theta_1$.
Note that  $\chi_1=0$ and $\left\{\frac{N}{M\hat i(\tilde c_{1, l})}\right\}=0$
by Theorem 3.7.
Let $\alpha_{1, i}=\frac{\phi_i}{\pi}$. Therefore we have 
\bea  \{m_1\alpha_{1, i}\}=\{A_{1, i}(N)\}=\left\{\left\{\frac{N\alpha_{1, i}}{\hat i(\tilde c_{1, l})}\right\}-\chi_{1, i}\right\}.\lb{4.72}
 \eea
 Hence 
 \bea\left\{\matrix{\{m_1\alpha_{1, i}\}<\epsilon\quad&&{\rm if}\quad (\chi_1,\,\chi_{1, i})=(0,\, 0), \cr
\{m_1\alpha_{1, i}\}>1-\epsilon\quad&&{\rm if}\quad (\chi_1,\,\chi_{1, i})=(0,\, 1), \cr}\right. \lb{4.73}\eea
 As in Case 1.1, we have by (\ref{3.26}), (\ref{4.43})   and Theorem 3.5
 \bea  &&i(\tilde c_{1, l})-2=2N+i(\tilde c_{1, l})-(2N+(3-1)\nn\\
= &&i(\tilde c_{1, l}^{2m_1+1})-i(\tilde c_{1, l}^{2m_1})\nn\\
=&&i(\tilde c_{1, l})-2+\sum_{1\le i\le 2}\left(2E\left(\frac{(2m_1+1)\phi_i}{2\pi}\right)-2E\left(\frac{2m_1\phi_i}{2\pi}\right)\right).\lb{4.74}\eea 
 This implies $E\left(\frac{(2m_1+1)\phi_i}{2\pi}\right)=E\left(\frac{2m_1\phi_i}{2\pi}\right) $
 for $i=1, 2$, and then  $\{m_1\alpha_{1, i}\}<\epsilon$ for $i=1, 2$ must hold. 
As in Case 1.1,  we can choose another $N^\prime\in\N$  with $m_i^\prime $ as defined in (\ref{4.70}) such that
 $ \{m_1^\prime\alpha_{1, i}\}>1-\epsilon$ for $i=1, 2$.
By (\ref{4.74}), we have 
\bea i(\tilde c_{1, l}^{2m_1^\prime})=i(\tilde c_{1, l}^{2m_1^\prime+1})
-(i(\tilde c_{1, l})-2+4)
= 2N^\prime+ i(\tilde c_{1, l})-i(\tilde c_{1, l})-2=2N^\prime-2.\lb{4.75}\eea
Hence  the closed geodesic found in Lemma 4.1 for $N^\prime$ must be $\tilde c_{2, l}$,
and then Step 1 holds in this case by Lemma 4.3.

{\bf Case 1.3.}  {\it We have $\Delta\notin\Q$ and $\hat i(\tilde c_{1, l})\notin\Q$.} 

In this case we have $P_{\tilde c_{1, l}}=R(\phi_1)\diamond R(\phi_2)$
with $\frac{\phi_i}{\pi}\notin\Q$ for $i=1, 2$ and $\phi_1=\theta_1$.
By Theorem 3.5, 
 we have  $T(\tilde c_{1, l})=1$ and then 
 \bea\hat \chi(\tilde c_{1, l}) =\chi(\tilde c_{1, l})= k_{\nu({\tilde c_{1, l}})}({\tilde c_{1, l}})\neq 0\lb{4.76}\eea
 by Proposition 2.1.
 
By Theorem 2.5  we have
\be n_1 \frac{\hat\chi(\tilde c_{1, 1})}{\hat{i}(\tilde c_{1, 1})}+n_2\frac{\hat\chi(\tilde c_{2, 1})}{\hat{i}(\tilde c_{2, 1})}=\sum_{1\le j\le 2, 1\le l\le n_j}\frac{\hat\chi(\tilde c_{j, l})}{\hat{i}(\tilde c_{j, l})}=B(3,1)=1.\lb{4.77}\ee
Thus we have   $\hat i(\tilde c_{2, l})\notin\Q$
also by $\hat i(\tilde c_{1, l})\notin\Q$ and (\ref{4.77}), and then Step 1 is true in his case
by Theorem 3.5.

{\bf Step 2.} {\it We have 
\bea M_k=b_k=\left\{\matrix{
    1,&\quad {\it if}\quad q=2  \cr
    2,&\quad {\it if}\quad q\in2N+2,  \cr
    0 &\quad {\it otherwise}. \cr}\right. \lb{4.78}\eea}

By Theorem 2.7, it is sufficient to show that $M_{2k+1}=0$ for $\forall k\in\Z$, i.e.,  
\bea \ol C_{2k+1}(E, \tilde c_{j, l}^m)=0,\qquad \forall  k\in\Z,\; 1\le j\le 2,\;1\le l\le n_j.\lb{4.79}
\eea
By Step 1, we can write
$P_{\tilde c_{j, l}}=R(\theta_j)\diamond M_j^\prime$ with  $\frac{\theta_j}{\pi}\notin\Q$ and $M_j^\prime\in \Sp(2)$ for $j=1, 2$.

By Lemma 4.1 and (\ref{4.53}), we have 
\bea \ol C_{q}(E, c_{j, l}^{2m_j})\neq 0,\quad q=2N\;{\rm or}\;2N+2.\lb{4.80}\eea

Note that by Theorem 3.5, the matrix $I_2$ and $N_1(1, 1)$ have the same
index iteration formula (\ref{3.20}) and can be viewed as $R(\theta)$ with $\theta=2\pi$.
although $I_2$ and $N_1(1, 1)$ have different nullity iteration formula (\ref{3.21}),
the discussion below for $I_2$ works also for  $N_1(1, 1)$. Thus in order to shorten the
length of the paper, we will not discuss $N_1(1, 1)$ separately.
Similarly the matrix $-I_2$ and $N_1(-1, 1)$ have the same
index iteration formula (\ref{3.20}) and can be viewed as $R(\theta)$ with $\theta=\pi$.
With this point of view, we have the following classification by Theorem 3.5:

{\bf Case 2.1.} {\it  $P_{\tilde c_{j, l}}=R(\theta_j)\diamond R(\varphi_j)$ with  $\frac{\varphi_j}{\pi}\in(0, 2]$.} 

By Theorem 3.5, we have 
\bea i(\tilde c_{j, l}^ m)&=&m(i(\tilde c_{j, l})-2)+
2E\left(\frac{m\theta_j}{2\pi}\right)
+2E\left(\frac{m\varphi_j}{2\pi}\right)-2,\nn\\
\nu(\tilde c_{j, l}^ m)&=&2\left(1-\varphi\left(\frac{m\varphi_j}{2\pi}\right)\right),\qquad\forall m\in\N,\lb{4.81}
\eea
and $i(\tilde c_{j, l})\in 2\N_0$.  

If $\frac{\varphi_j}{2\pi}\notin\Q$,  we have  $i(\tilde c_{j, l}^ m)\in 2\N_0$ and $\nu(\tilde c_{j, l}^ m)=0$ for $m\in\N$. Hence  (\ref{4.79}) holds by Proposition 2.1.

If $\frac{\varphi_j}{2\pi}\in\Q$, write $\frac{\varphi_j}{2\pi}=\frac{r}{s}$ with $r, s\in\N$ and $(r, s)=1$. Then we have  $i(\tilde c_{j, l}^ m)$ is always even and $\nu(\tilde c_{j, l}^ m)=2$
if $s|m$,  $\nu(\tilde c_{j, l}^ m)=0$ otherwise.
Thus by (\ref{4.80}) and Proposition 2.1,   We have $k_0(\tilde c_{j, l}^ {2m_j})\neq 0$
or $k_2(\tilde c_{j, l}^ {2m_j})\neq 0$, and then $k_1(\tilde c_{j, l}^ {2m_j})=0$
by  Proposition 2.2. By Proposition 2.3, we have
\bea k_l(\tilde c_{j, l}^ m)=k_l(\tilde c_{j, l}^ {2m_j}),\quad s|m,\;\forall l\in\Z.
\lb{4.82}\eea
Hence  (\ref{4.79}) holds by Proposition 2.1 and (\ref{4.82}).

{\bf Case 2.2.} {\it  $P_{\tilde c_{j, l}}=R(\theta_j)\diamond N_1(1, -1)$.}
By Theorem 3.5, we have 
\bea i(\tilde c_{j, l}^ m)&=&m(i(\tilde c_{j, l})-1)+
2E\left(\frac{m\theta_j}{2\pi}\right)
-1,\nn\\
\nu(\tilde c_{j, l}^ m)&=&1,\qquad\forall m\in\N,\lb{4.83}
\eea
and $i(\tilde c_{j, l})\in 2\N-1$.  

Then we have  $i(\tilde c_{j, l}^ m)$ is always odd.
Thus by (\ref{4.80}) and Proposition 2.1,   We have $k_1(\tilde c_{j, l}^ {2m_j})\neq 0$
and then $k_0(\tilde c_{j, l}^ {2m_j})=0$
by  Proposition 2.2. By Proposition 2.3, we have
\bea k_l(\tilde c_{j, l}^ m)=k_l(\tilde c_{j, l}^ {2m_j}),\quad m\in\N,\;\forall l\in\Z.
\lb{4.84}\eea
Hence  (\ref{4.79}) holds by Proposition 2.1 and (\ref{4.84}).

{\bf Case 2.3.} {\it  $P_{\tilde c_{j, l}}=R(\theta_j)\diamond N_1(-1, 1)$.}
By Theorem 3.5, we have 
\bea i(\tilde c_{j, l}^ m)&=&m(i(\tilde c_{j, l})-1)+
2E\left(\frac{m\theta_j}{2\pi}\right)
-1,\nn\\
\nu(\tilde c_{j, l}^ m)&=&\frac{1+(-1)^m}{2}\qquad\forall m\in\N,\lb{4.85}
\eea
and $i(\tilde c_{j, l})\in 2\N_0$.  

Then we have  $i(\tilde c_{j, l}^ m)\in2\N_0$ and $\nu(\tilde c_{j, l}^ m)=0$
if $m\in2\N-1$;  $i(\tilde c_{j, l}^ m)\in2\N-1$ and $\nu(\tilde c_{j, l}^ m)=1$
if $m\in2\N$.  
Thus by (\ref{4.80}) and Proposition 2.1,   We have $k_1(\tilde c_{j, l}^ {2m_j})\neq 0$
and then $k_0(\tilde c_{j, l}^ {2m_j})=0$
by  Proposition 2.2. By Proposition 2.3, we have
\bea k_l(\tilde c_{j, l}^ m)=k_l(\tilde c_{j, l}^ {2m_j}),\quad m\in2\N,\;\forall l\in\Z.
\lb{4.86}\eea
Hence  (\ref{4.79}) holds by Proposition 2.1 and (\ref{4.86}).

{\bf Case 2.4.} {\it  $P_{\tilde c_{j, l}}=R(\theta_j)\diamond H$ with $H$ being hyperbolic.}

By Theorem 3.5, we have 
\bea i(\tilde c_{j, l}^ m)&=&m(i(\tilde c_{J, l})-1)+
2E\left(\frac{m\theta_j}{2\pi}\right)
-1,\nn\\
\nu(\tilde c_{j, l}^ m)&=&0,\lb{4.87}
\eea
and $i(\tilde c_{j, l})\in \N_0$.  

If $i(\tilde c_{j, l})\in 2\N-1$, then we have  $i(\tilde c_{j, l}^ m)\in2\N-1$ 
for $m\in\N$. This contradicts to (\ref{4.80}) by Proposition 2.1. Hence this case
can not appear. 

If $i(\tilde c_{j, l})\in 2\N_0$, then we have  $i(\tilde c_{j, l}^ m)\in2\N_0$ if $m\in2\N-1$; $i(\tilde c_{j, l}^ m)\in2\N-1$ if $m\in2\N$.  This contradicts to (\ref{4.80}) by Proposition 2.1. Hence this case
can not appear also.

In particular, we have $M_2=b_2=1$, therefore there must be a prime closed geodesic
$\tilde c_{\alpha, l}$ on  $(\tilde M, \tilde F)$ such that 
\bea i(\tilde c_{\alpha, l})=2, \quad k_0(\tilde c_{\alpha, l})=1,\lb{4.88}\eea
for some $\alpha\in\{1, 2\}$. In fact, $M_2=1$ implies there exist
some  $\alpha\in\{1, 2\}$ and $m\in\N$ such that $\ol C_2(E, \tilde c_{\alpha, l}^m)=1$.
Thus by (\ref{4.2}) and Proposition 2.1, we have i$(\tilde c_{\alpha, l}^m)=2$ and $k_0(\tilde c_{\alpha, l}^m)=1$. Now if $m=1$, then (\ref{4.88}) is true. If $m>1$, then  
$i(\tilde c_{\alpha, l}^m)\ge i(\tilde c_{\alpha, l}^{m-1})+\nu(\tilde c_{\alpha, l}^{m-1})
\ge\cdots\ge i(\tilde c_{\alpha, l})+\nu(\tilde c_{\alpha, l})>2$ by (\ref{4.2}) and (\ref{4.3}) provided 
$\nu(\tilde c_{\alpha, l})>0$,  this contradiction implies $\nu(\tilde c_{\alpha, l})=0$, 
and then (\ref{4.88}) holds by Proposition 2.1.

Thus we have 
\bea \ol C_{2N+2}(E, \tilde c_{\alpha, l}^{2m_\alpha+1})=\Q,
\lb{4.89}\eea
by (\ref{3.26}), (\ref{4.88}) and Proposition 2.3.

By Lemma 4.1, we have 
\bea i(\tilde c_{1, l}^{2m_{1}})+\nu(\tilde c_{1, l}^{2m_{1}})=2N+2,\quad
\ol{C}_{2N+2}(E, \tilde c_{1, l}^{2m_{1}})\neq 0,\lb{4.90}\eea
 for $1\le l\le n_{1}$. 
 Thus we have 
 \bea 2=b_{2N+2}=M_{2N+2}\ge 1+n_1.\lb{4.91}
\eea
This implies that $n_1=1$.

{\bf Step 3.}  {\it We have $n_2=1$ provided $\tilde c_{1, l}$ belongs to 
Cases 1.1 or 1.2 in Step 1.}

In fact, by  Cases 1.1 and 1.2, we can find $(N^\prime, m_1^\prime, m_2^\prime)$ such that
\bea i(\tilde c_{2, l}^{2m_2^\prime})+\nu(\tilde c_{2, l}^{2m_2^\prime})=2N^\prime+2,\quad
\ol{C}_{2N^\prime+2}(E, \tilde c_{2, l}^{2m_2^\prime})\neq 0,\lb{4.92}\eea
 for $1\le l\le n_2$.  As above, we have $\ol C_{2N^\prime+2}(E, \tilde c_{\alpha, l}^{2m^\prime_\alpha+1})=\Q$. Then as in (\ref{4.91}), we have $2=b_{2N^\prime+2}=M_{2N^\prime+2}\ge 1+n_2$. This implies $n_2=1$.

{\bf Step 4.}  {\it We have $n_2=1$ provided $\tilde c_{1, l}$ belongs to 
Case 1.3 in Step 1.}

Firstly we show
\bea i(\tilde c_{1, l}^{2m_1-m})< 2N-2,\quad \forall m\ge 2.\lb{4.93}\eea

In fact, we have $P_{\tilde c_{1, l}}=R(\phi_1)\diamond R(\phi_2)$
with $\frac{\phi_i}{\pi}\notin\Q$ for $i=1, 2$.
 By (\ref{4.3}) and (\ref{4.6}), if $i(\tilde c_{1, l})>2$,
then (\ref{4.93}) holds. Thus it remains to consider the case
$i(\tilde c_{1, l})=2$. By Lemma 2 of \cite{Rad5} and Theorem 3.5, we have
\bea \hat i(\tilde c_{1, l})&=&i(\tilde c_{1, l})+p_-+p_0-r+\sum_{i=1}^r
\frac{\phi_i}{\pi}\nn\\
&=&i(\tilde c_{1, l})-2+\sum_{i=1}^2\frac{\phi_i}{\pi}
>2.\lb{4.94}\eea
Plugging
$i(\tilde c_{1, l})=2$ into (\ref{4.94}) yields
\bea  \sum_{i=1}^2\left(\frac{\phi_i}{\pi}-1\right)>0.\lb{4.95}\eea
Thus we may assume without loss of generality that $\phi_1\in (\pi, 2\pi)$.
By Theorem 3.5 we have
\bea i(\tilde c_{1, l}^m)=m(i(\tilde c_{1, l})-2)+2\sum_{i=1}^2{E}\left(\frac{m\phi_i}{2\pi}\right)-2
=2\sum_{i=1}^2{E}\left(\frac{m\phi_i}{2\pi}\right)-2.
\lb{4.96}\eea
By (\ref{4.3}) and (\ref{4.6}), in order to prove (\ref{4.93}), it is sufficient
to prove
\bea i(\tilde c_{1, l}^{2m_1-2})<i(\tilde c_{1, l}^{2m_1-1}).\lb{4.97}\eea
By (\ref{4.96}), in order to prove (\ref{4.97}), it is sufficient
to prove
\bea {E}\left(\frac{(2m_1-2)\phi_1}{2\pi}\right)
<{E}\left(\frac{(2m_1-1)\phi_1}{2\pi}\right).
\lb{4.98}\eea
In order to satisfy (\ref{4.98}), it is sufficient
to choose
\bea \delta<\min\left\{\frac{\phi_1}{\pi}-1,\,1-\frac{\phi_1}{2\pi}\right\},
\lb{4.99}\eea
where $\delta$ is given by (\ref{3.29}). This proves  (\ref{4.93}).

Thus by (\ref{4.11}), (\ref{4.14}), (\ref{4.80}), (\ref{4.93}) and Proposition 2.1, we have 
\bea \ol{C}_{2N-2}(E, \tilde c_{1, l}^{2m_1-1})=\Q,
\quad \ol{C}_{2N-2}(E, \tilde c_{1, l}^{m})=0,\;\forall m\neq 2m_1-1.\lb{4.100}\eea
Thus we have 
 \bea 2=b_{2N-2}=M_{2N-2}\ge 1+n_2.\lb{4.101}
\eea
This implies that $n_2=1$.

Now we can complete the proof of Theorem 1.5 as follows:
by Steps 2-4 and Proposition 2.1, we have 
 \bea 2=b_{2N}=M_{2N}=1.\lb{4.102}
\eea
In fact, by the proof of Step 2, we have $\ol C_q(E, \tilde c_{j,  l}^{2m_j})=\Q$ for 
$q=2N$ or $2N+2$ by Proposition 2.2, thus (\ref{4.102}) holds.
This contradiction proves Theorem 1.5.
\hfill\hb

\bibliographystyle{abbrv}

\bigskip

\end{document}